\documentclass[11pt, reqno]{amsart}   	
\usepackage{geometry}                		
\usepackage{graphicx}				
\usepackage{amssymb}
\usepackage{amsmath}
\usepackage{comment}
\usepackage{amsthm}
\usepackage{mathrsfs}
\allowdisplaybreaks
\newtheorem{theorem}{Theorem}
\numberwithin{theorem}{section}
\newtheorem{proposition}[theorem]{Proposition}
\newtheorem{lemma}[theorem]{Lemma}
\newtheorem{corollary}[theorem]{Corollary}
\newcommand{\EndProof}{\hfill $\square$}
\usepackage[colorlinks=true, linkcolor=blue,
 citecolor=blue, urlcolor=blue]{hyperref}
 \usepackage{dsfont}
 \usepackage{tikz-cd}
 \usepackage{caption}

\date{}							

\begin{document}
\title{A Hecke Algebra on the Double Cover of A Chevalley Group Over $\mathbb{Q}_{2}$} 
\author{Edmund Karasiewicz}
\email{karasiew@post.bgu.ac.il}
\address{Edmund Karasiewicz: Department of Mathematics, Ben Gurion University of the Negev, Be'er Sheva,  Israel 8410501}
\begin{abstract} We prove that a certain genuine Hecke algebra $\mathcal{H}$ on the non-linear double cover of a simple, simply-laced, simply-connected, Chevalley group $G$ over $\mathbb{Q}_{2}$ admits a Bernstein presentation. This presentation has two consequences. First, the Bernstein component containing the genuine unramified principal series is equivalent to $\mathcal{H}$-mod. Second, $\mathcal{H}$ is isomorphic to the Iwahori-Hecke algebra of the linear group $G/Z_{2}$, where $Z_{2}$ is the $2$-torsion of the center of $G$. This isomorphism of Hecke algebras provides a correspondence between certain genuine unramified principal series of the double cover of $G$ and the Iwahori-unramified representations of the group $G/Z_{2}$.
\end{abstract}
\maketitle

\section{Introduction} 
Let $\textbf{G}$ be a simple, simply-laced, simply-connected, Chevalley group over a $p$-adic field $F$, where $F$ contains the $n$-th roots of unity $\mu_{n}$. Let $\textbf{T}\subset \textbf{B}\subset \textbf{G}$ be a maximal torus and Borel subgroup, respectively. If $\textbf{H}$ is an algebraic group, we will write $H$ for the $F$-points of $\textbf{H}$. Let $\widetilde{G}$ be the $n$-fold cover of $G$. The covering group $\widetilde{G}$ fits into a central extension
\begin{equation}\label{CentralExt1}
1\rightarrow \mu_{n} \rightarrow \widetilde{G}\stackrel{\mathrm{pr}}{\rightarrow} G\rightarrow 1.
\end{equation}
We fix an embedding $\varepsilon:\mu_{n}\rightarrow \mathbb{C}^{\times}$. 

When $\mathrm{GCD}(n,p)=1$, the tame case, Savin \cite{S88,S04} studied the genuine Iwahori-Hecke algebra of $\widetilde{G}$ and its relation to the representation theory of $\widetilde{G}$. Among his results is a Bernstein presentation for the $\varepsilon$-genuine Iwahori-Hecke algebra, which yields the following two consequences. First, the Bernstein component of $\varepsilon$-genuine unramified principal series is equivalent to the category of modules of the $\varepsilon$-genuine Iwahori-Hecke algebra. Second, there is an equivalence of categories between the Bernstein component of $\varepsilon$-genuine Iwahori-unramified representations of $\widetilde{G}$ and the Bernstein component of Iwahori-unramified representations of $G^{\prime}=G/Z_{n}$, where $Z_{n}$ is the $n$-torsion in the center of $G$. (For details on Bernstein components see \cite{BD84,BK98}.)

When GCD$(n,p)\neq1$, the wild case, less is known. Loke-Savin \cite{LS10b} established the analogous results when $G=\mathrm{SL}(2,\mathbb{Q}_{2})$ and $n=2$. Wood \cite{W14} rephrased these results in terms of the even Weil representation and extended them to the case of $G=$Sp$(2r,\mathbb{Q}_{2})$ and $n=2$. Takeda-Wood \cite{TW18} reproved the results of Wood for any 2-adic field, and proved an analogous result for the odd Weil representation.

At this point, we describe how the GCD$(n,p)$ affects matters. If GCD$(n,p)=1$, then $\mathrm{pr}^{-1}(I)\cong I\times\mu_{n}$, where $I$ is an Iwahori subgroup of $G$. Thus $I$ can be embedded as a subgroup in $\widetilde{G}$. We will say that $I$ \textit{splits} sequence (\ref{CentralExt1}) and call the embedding a \textit{splitting} of $I$ into $\widetilde{G}$. Consequently, $\widetilde{G}$ possesses an $\varepsilon$-genuine Iwahori-Hecke algebra. When GCD$(n,p)\neq 1$, $I$ need not split sequence (\ref{CentralExt1}), leaving no obvious candidate for the analog of the $\varepsilon$-genuine Iwahori-Hecke algebra of $\widetilde{G}$. However, Loke-Savin \cite{LS10b} identify a suitable alternative when $G=SL(2,\mathbb{Q}_{2})$ and $n=2$. 

In this paper we extend the ideas of Loke-Savin \cite{LS10b} to the double cover of higher-rank simple, simply-laced, simply-connected, Chevalley groups over $\mathbb{Q}_{2}$. In particular, $n=2$ and $F=\mathbb{Q}_{2}$. (Since $n=2$, $\varepsilon$ is unique and will be omitted.) 

We begin by identifying a Hecke algebra $\mathcal{H}$ to replace the Iwahori-Hecke algebra. The definition of $\mathcal{H}$ involves the following compact open subgroups and representation. Let $\Gamma_{0}(4)$ be the pre-image of $\textbf{B}(\mathbb{Z}/4\mathbb{Z})$ under the mod $4$ reduction map $\textbf{G}(\mathbb{Z}_{2})\rightarrow \textbf{G}(\mathbb{Z}/4\mathbb{Z})$. Let $\Gamma_{1}(4)$ be the pre-image of the unipotent radical of $\textbf{B}(\mathbb{Z}/4\mathbb{Z})$. It is crucial that $\Gamma_{1}(4)$ splits sequence (\ref{CentralExt}) (Theorem \ref{Gamma14Splitting}). Fix $(\tau,E)$, an irreducible genuine $\mathrm{pr}^{-1}(\Gamma_{0}(4))$-representation that is trivial on the image of the splitting of $\Gamma_{1}(4)$. The principal Hecke algebra of the present paper is
\[\mathcal{H}\stackrel{\mathrm{def}}{=}\{f\in C^{\infty}_{c}(\widetilde{G},\text{End}(E))| f( k_1gk_2)=\tau(k_1)f(g)\tau(k_2),\,\text{for all }k_1,k_2\in \mathrm{pr}^{-1}(\Gamma_0(4))\}.\]

In Theorem \ref{HeckeIso},  we prove that $\mathcal{H}$ admits a Bernstein presentation. As in the tame case, this presentation yields two consequences. 

\begin{theorem}(See Theorem \ref{BernsteinBlock}.)\label{BernsteinBlockIntro} Let $(\pi,\mathscr{V})$ be a smooth $\widetilde{G}$-representation that is generated by its $\tau$-isotypic vectors
\begin{enumerate}
\item\label{BB1Intro}  Every subquotient of $(\pi,\mathscr{V})$ is generated by its $\tau$-isotypic vectors.

\item\label{BB2Intro} If $(\pi,\mathscr{V})$ is irreducible, then there is an unramified character $\chi:T\rightarrow \mathbb{C}^{\times}$ such that $(\pi,\mathscr{V})$ is isomorphic to an irreducible subquotient of $\mathrm{Ind}_{\widetilde{B}}^{\widetilde{G}}(i(\chi))$. (The definition of $i(\chi)$ can be found in Subsection \ref{CoveringTorus}.)
\end{enumerate}
\end{theorem}

Let $I^{\prime}$ be an Iwahori subgroup of $G^{\prime}=G/Z_{2}$.

\begin{theorem}(See Theorem \ref{main}.)
The Bernstein presentation of $\mathcal{H}$ induces an algebra isomorphism between $\mathcal{H}$ and the Iwahori-Hecke algebra of $G^{\prime}$. This isomorphism induces an equivalence of categories between the category of smooth $G^{\prime}$-representations generated by their $I^{\prime}$-fixed vectors and the category of smooth $\widetilde{G}$-representations generated by their $\tau$-isotypic vectors.
\end{theorem}
Our approach to achieve these results is indebted to Savin \cite{S04} and Loke-Savin \cite{LS10a, LS10b}. 

Now we describe the contents of this paper in more detail. In Section \ref{notation} we establish notation. 

Section \ref{splitting} contains a proof of the splitting of $\Gamma_{1}(4)$ (Theorem \ref{Gamma14Splitting}). 
In addition to existence, we prove an important technical result (Proposition \ref{SplitWeylConjugate}) required to show that certain functions, which constitute a basis of the Hecke algebra, are well-defined.

In Section \ref{doublecosets} we isolate some double coset calculations, which provide a point of departure for our study of $\mathcal{H}$. The main result in this section (Proposition \ref{CosetUniRed}) describes a preferred representative of each double coset.

Section \ref{LinearHeckeAlg} introduces $\underline{\mathcal{H}}$ a Hecke algebra on the linear group $G$. This linear Hecke algebra is significant for two resasons. First, multiplication in $\mathcal{H}$ is related to multiplication in $\underline{\mathcal{H}}$, where the existence of a particular algebra homomorphism can simplify computations. This can be seen in Proposition \ref{LengthMult}. Second, we can use $\underline{\mathcal{H}}$ to prove Theorem \ref{EqualLengths}. This theorem, which was pointed out to me by Gordan Savin, provides a starting point for showing that $\mathcal{H}$ has an Iwahori-Matsumoto presentation (Theorem \ref{AffineHeckePresentation}).

Section \ref{2adicHecke}, the principal section of this paper, contains our study of the Hecke algebra $\mathcal{H}$, which culminates in a Bernstein presentation of $\mathcal{H}$ (Theorem \ref{HeckeIso}). Now we will briefly outline the argument.

First, We begin with support calculations (Proposition \ref{NoUni}, Proposition \ref{ModCoCharRed}), which ultimately allow us to construct a $\mathbb{C}$-basis of the Hecke algebra (Proposition \ref{BasisWellDefined}). This step breaks up into two pieces. First, we must show that certain double cosets do not support any functions in the Hecke algebra. Second, we construct a basis for $\mathcal{H}$, invoking Proposition \ref{SplitWeylConjugate} to prove that the basis functions are well defined.

Second, we show that $\mathcal{H}$ satisfies the braid relations (Proposition \ref{LengthMult}) and quadratic relations (Proposition \ref{QuadRel}). 
Using these relations and standard facts about affine Hecke algebras we prove that $\mathcal{H}$ admits an Iwahori-Matsumoto presentation (Proposition \ref{AffineHeckePresentation}). 

Third, we apply the results of Lusztig \cite{L89} to prove that  $\mathcal{H}$ satisfies the Bernstein relations, and use Savin's trick (Lemma 7.6, \cite{S04}) to show that these relations imply all others. This completes the proof of Theorem \ref{HeckeIso}.  

Section \ref{BernsteinComponent} contains a proof that the Bernstein component containing the genuine unramified principal series is equivalent to $\mathcal{H}$-mod (Theorem \ref{BernsteinBlock}). 

In Section \ref{ShimuraCor}, we describe the resulting local Shimura correspondence. The presentation proved in Theorem \ref{HeckeIso} agrees with the Bernstein presentation of an affine Hecke algebra. This leads to an isomorphism between $\mathcal{H}$ and the Iwahori-Hecke algebra of $G/Z_{2}$, which yields the desired Shimura correspondence (Theorem \ref{main}).

The construction of the isomorphism between $\mathcal{H}$ and the Iwahori-Hecke Algebra depends on several choices. We conclude the present work by enumerating these choices.\\


\noindent\textbf{Acknowledgements:} I would like to thank Martin Weissman for suggesting this project and his guidance throughout the preparation of this article. I would also like to thank Gordan Savin for his extensive comments on an earlier draft of this work. In particular, he brought Theorem \ref{EqualLengths} to my attention, which greatly simplified the proof of Theorem \ref{HeckeIso}.

\section{Notation}\label{notation}

\subsection{Root System}

Let $(\Phi,\Delta,\mathcal{E})$ be a reduced irreducible simply-laced root system in the Euclidean vector space $\mathcal{E}$ with roots $\Phi$ and simple roots $\Delta$. Let $r=|\Delta|$. We write $\Phi^{+}$ for the set of positive roots and $\Phi^{-}=-\Phi^{+}$ for the negative roots. Let $\varrho=\frac{1}{2}\sum_{\alpha\in \Phi^{+}}\alpha$. For $\alpha, \beta\in\Phi$ we say that $\alpha\succ\beta$ if $\alpha-\beta$ is a sum of positive roots. Let $\alpha_{1}\ldots,\alpha_{r}$ be an enumeration of the simple roots and let $\alpha_{0}$ be the lowest root in $\Phi$. Associated with this root system there is a semi-simple complex lie algebra $\mathfrak{g}$.

\subsection{Chevalley Group}

By choosing a Chevalley basis for $\mathfrak{g}$ we can construct $\textbf{G}$, the associated simply connected Chevalley group over $\mathbb{\mathbb{Z}}$ with maximal torus $\textbf{T}$. Let $\textbf{B}\supset \textbf{T}$ denote the Borel subgroup associated to $\Delta$ with unipotent radical $\textbf{U}$; let $\textbf{B}_{-}$ denote the Borel subgroup opposite to $\textbf{B}$ with unipotent radical $\textbf{U}_{-}$. Let $\textbf{N}$ denote the normalizer of $\textbf{T}$ in $\textbf{G}$, and let $W$ be the Weyl group of $\textbf{G}$ with respect to $\textbf{T}$.

The torus $\textbf{T}$ has a group of rational characters $X=X^*(\textbf{T})$ and cocharacters $Y=X_*(\textbf{T})$. Let $\langle\cdot,\cdot\rangle:X\times X \rightarrow \mathbb{Z}$ be the Killing form normalized so that the roots have length $2$. Using the Killing Form we can identify $X$ and $Y$ and the root $\alpha$ is identified with the coroot $\alpha^{\vee}$ (since our root system is simply-laced). This identification induces a bilinear form on $Y$, also denoted by $\langle\cdot,\cdot\rangle$. We define the modified cocharacter lattice to be $\widetilde{Y}=\{y\in Y|\langle y,y^\prime\rangle\in2\mathbb{Z}\text{ for all }y^\prime\in Y\}$. If $Y^*$ is the lattice dual to $Y$ in $Y\otimes \mathbb{R}$ with respect to the Killing form, then $\widetilde{Y}=Y\cap2Y^*$.

Let $G=\textbf{G}(\mathbb{Q}_2)$, the $\mathbb{Q}_2$-rational points of $\textbf{G}$. Similarly, we will write $T,U,U_{-}, B,B_{-},\ldots$ for the $\mathbb{Q}_{2}$-points of $\textbf{T},\textbf{U},\textbf{U}_{-},\textbf{B},\textbf{B}_{-},\ldots$, respectively. 

We will be interested in the topological two-fold cover of $G$, but first we recall some facts about the universal central extension of $G$.  Let $\mathrm{St}(\Phi,\mathbb{Q}_{2})$ denote the universal central extension of $G$ over $\mathbb{Q}_{2}$. The group $\mathrm{St}(\Phi,\mathbb{Q}_{2})$ is generated by the elements $x_{\alpha}^{\prime}(t)$, where $\alpha\in\Phi$ and $t\in\mathbb{Q}_{2}$, and subject to the relations 

\begin{equation*}
\begin{array}{ccl}
x^{\prime}_{\alpha}(t)x^{\prime}_{\alpha}(u)&=&x_{\alpha}^{\prime}(t+u),\\
\lbrack x_{\alpha}^{\prime}(t),x_{\beta}^{\prime}(u)\rbrack&\stackrel{\mathrm{def}}{=}&x_{\alpha}^{\prime}(t)x_{\beta}^{\prime}(u)x_{\alpha}^{\prime}(-t)x_{\beta}^{\prime}(-u)=
\begin{cases}
1, \text{ if }\alpha+\beta\notin \Phi;\\
x_{\alpha+\beta}^{\prime}(c(\alpha,\beta) tu), \text{ if }\alpha+\beta\in\Phi,
\end{cases}
\end{array}
\end{equation*}
where $c(\alpha,\beta)=\pm1$. (Recall that $G$ is simply-laced.) For more details see \cite{S68}, where Steinberg writes $G^{\prime}$ instead of $\mathrm{St}(\Phi,\mathbb{Q}_{2})$.

Furthermore, for $\alpha\in\Phi$ and $t\in\mathbb{Q}_{2}^{\times}$ we let 
\begin{equation*}
\begin{array}{rcl}
w^{\prime}_{\alpha}(t)&=&x^{\prime}_{\alpha}(t)x^{\prime}_{-\alpha}(-t^{-1})x^{\prime}_{\alpha}(t),\\ 
h^{\prime}_{\alpha}(t)&=&w^{\prime}_{\alpha}(t)w^{\prime}_{\alpha}(-1).
\end{array}
\end{equation*}
The work of Moore \cite{M68} and Matsumoto \cite{M69} provides a presentation for the kernel of the central extension
\begin{equation}\label{UCE}
1\rightarrow \mathrm{Ker}(\mathrm{pr}^{\prime})\rightarrow \mathrm{St}(\Phi,\mathbb{Q}_{2})\stackrel{\mathrm{pr}^{\prime}}{\rightarrow}G\rightarrow 1.
\end{equation}
The elements of the form $\{t,u\}\stackrel{\text{def}}{=}h^{\prime}_{\alpha}(t)h^{\prime}_{\alpha}(u)h^{\prime}_{\alpha}(tu)^{-1}$ generate $\mathrm{Ker}(\mathrm{pr}^{\prime})$ and satisfy the relations described in Theorem 12 in \cite{S68}, where Steinberg writes $f(t,u)$ for $\{t,u\}$.

By \cite{M68}, the push-out of sequence (\ref{UCE}) via the quadratic Hilbert symbol $(\cdot,\cdot)_{2}:\mathrm{Ker}(\mathrm{pr}^{\prime})\rightarrow \mu_{2}=\{\pm1\}$ yields the group $\widetilde{G}$, the unique nontrivial topological two-fold central extension of $G$.  In particular, we have the following commutative diagram with exact rows:
\begin{equation}
\begin{tikzcd}
1\arrow[r] &\mathrm{Ker}(\mathrm{pr}^{\prime})\arrow[d,"{(\cdot,\cdot)_{2}}"]\arrow[r] &\mathrm{St}(\Phi,\mathbb{Q}_{2})\arrow[d,"{\mathrm{pr}^{\prime\prime}}"] \arrow[r,"\mathrm{pr}^{\prime}"]&G\arrow[d]\arrow[r] &1\\
1\arrow[r] &\mu_2\arrow[r] &\widetilde{G}\arrow[r,"\mathrm{pr}"] &G\arrow[r] &1.
\end{tikzcd}
\end{equation}

For each $\alpha\in \Phi$ and $t\in\mathbb{Q}_{2}$ let $\tilde{x}_\alpha(t)=\mathrm{pr}^{\prime\prime}(x^{\prime}_{\alpha}(t))$. For $\alpha\in\Phi$ and $t\in\mathbb{Q}_{2}^{\times}$ define $\tilde{w}_{\alpha}(t)=\mathrm{pr}^{\prime\prime}(w^{\prime}_{\alpha}(t))$ and $\tilde{h}_{\alpha}(t)=\mathrm{pr}^{\prime\prime}(h^{\prime}_{\alpha}(t))$. Similarly, let $x_\alpha(t)=\mathrm{pr}(\tilde{x}_\alpha(t))$, $w_{\alpha}(t)=\mathrm{pr}(\tilde{w}_{\alpha}(t)),$ and $h_{\alpha}(t)=\mathrm{pr}(\tilde{h}_{\alpha}(t))$. 

For a ring $R$, let $\mathbf{U}^{*}(R)$ be the subgroup of $\widetilde{G}$ generated by the elements $\tilde{x}_{\alpha}(t)$, where $t\in R$ and $\alpha\in\Phi^{+}$; define $\mathbf{U}^{*}_{-}(R)$ similarly. Let $T^{*}_{1}$ be the subgroup of $\widetilde{G}$ generated by $\tilde{h}_{\alpha}(t)$, where $t\in 1+4\mathbb{Z}_{2}$ and $\alpha\in\Phi$, and let $T_{1}=\mathrm{pr}(\widetilde{T}_1)$. 

There are a few subgroups $J$ of $G$ which possess a splitting of the sequence 
\begin{equation}
\begin{tikzcd}
1\arrow[r] &\mu_2\arrow[r] &\widetilde{G}\arrow[r,"\mathrm{pr}"] &G\arrow[r] &1\label{CentralExt},
\end{tikzcd}
\end{equation}
in other words, a group homomorphism $f:J\rightarrow \widetilde{G}$ such that $\mathrm{pr}\circ f=\mathrm{id}_{H}$. The following maps split sequence (\ref{CentralExt}).
\begin{align}
\label{CanSpl1}\mathbf{U}(R)\rightarrow\,& \mathbf{U}^{*}(R), \text{ defined by } x_{\alpha}(t)\mapsto \tilde{x}_{\alpha}(t)\text{, for }\alpha\in\Phi^{+},\,t\in R;\\ 
\label{CanSpl2}\mathbf{U}_{-}(R)\rightarrow\,& \mathbf{U}^{*}_{-}(R), \text{ defined by } x_{\alpha}(t)\mapsto \tilde{x}_{\alpha}(t)\text{, for }\alpha\in\Phi^{-},\,t\in R;\\ 
\label{CanSpl3}T_{1}\rightarrow\,& T^{*}_{1},\text{ defined by } h_{\alpha}(t)\mapsto \tilde{h}_{\alpha}(t)\text{, for }\alpha\in\Phi,\,t\in 1+4\mathbb{Z}_{2}.
\end{align}

Note that the Steinberg relations and the fact that $(2,2)_{2}=1$ imply that the subgroup of $\widetilde{G}$ generated by $h_{\alpha}(2)$ for all $\alpha\in \Delta$ also splits the exact sequence (\ref{CentralExt}) and is isomorphic to $Y$. For $\lambda=\sum_{j} c_{j}\alpha_{j}\in Y$, let $2^\lambda$ denote $\prod_{j} h_{\alpha_{j}}(2)^{c_{j}}$. Let $\Upsilon:Y\rightarrow\widetilde{G}$ be the map defined by $\lambda\mapsto 2^{\lambda}$. Note that $\Upsilon$ is a group isomorphism. 

Let $\mathcal{W}$ be the subgroup of $\widetilde{G}$ generated by the elements $\tilde{w}_{\alpha}(1)$, where $\alpha\in \Phi$. Let $\widetilde{N}^{\prime}$ be the subgroup of $\widetilde{G}$ generated by the elements $\tilde{w}_{\alpha}(1)$ and $2^{\lambda}$, where $\alpha\in\Phi$ and $\lambda\in\widetilde{Y}$. Using the Steinberg relations one can show that $\widetilde{N}^{\prime}\cong\mathcal{W}\ltimes \widetilde{Y}$.

Consider the map $\textbf{G}(\mathbb{Z}_2)\rightarrow \textbf{G}(\mathbb{Z}_{2}/2^{k}\mathbb{Z}_{2})$ defined by reduction modulo $2^{k}$. Let $\Gamma(2^k)$ be the kernel of this map. Let $\Gamma_{0}(2^k)$ be the inverse image of $\textbf{B}(\mathbb{Z}/2^k\mathbb{Z})$ and let $\Gamma_{1}(2^k)$ to be the inverse image of $\textbf{U}(\mathbb{Z}/2^k\mathbb{Z})$. 

Given a subgroup $J\subseteq G$, let $\widetilde{J}=\mathrm{pr}^{-1}(J)$. A representation of $\widetilde{J}$ is said to be genuine if $\mu_{2}\subset \widetilde{J}$ acts nontrivially.

\subsection{Affine Weyl Group} 

Two affine Weyl groups are pertinent to our study. The first is $W_{\mathrm{aff}}$, associated to the root system $(\Phi,\mathcal{E})$; the second will be an extended affine Weyl group $\widetilde{W}_{\mathrm{aff}}$ associated to $(\frac{1}{2}\Phi,\mathcal{E})$. We begin with $W_{\mathrm{aff}}$. 

Given a root system $(\Phi,\mathcal{E})$, there is an associated affine Weyl group $W_{\mathrm{aff}}$ generated by the reflections through the hyperplanes $H_{\alpha,k}\stackrel{\mathrm{def}}{=}\{v\in\mathcal{E}|\langle\alpha,v\rangle = k\}$, where $\alpha\in \Phi$ and $k\in \mathbb{Z}$. Let $w_{\alpha,k}$ be the reflection fixing $H_{\alpha,k}$ defined by $w_{\alpha,k}(v)=v-\langle\alpha,v\rangle\alpha^{\vee}$. The reflections $w_{\alpha_{0},-1},w_{\alpha_{1},0}\ldots,w_{\alpha_{r},0}$ are a set of Coxeter generators for $W_{\mathrm{aff}}$. We will write $w_{\alpha_{i}}$ for $w_{\alpha_{i},0}$, where $i=1,\ldots, r$, and $w_{\alpha_{0}}^{\prime}$ for $w_{\alpha_{0},-1}$.

The group $W_{\mathrm{aff}}$ can be decomposed as a semi-direct product $W_{\mathrm{aff}}\cong  W\ltimes Y$ (since $Y$ is the lattice of coroots), and it can be realized as the quotient $\textbf{N}(\mathbb{Q}_{2})/\textbf{T}(\mathbb{Z}_{2})$. Furthermore, there is a section $\textbf{s}:W_{\mathrm{aff}}\cong W\times Y \rightarrow \textbf{N}(\mathbb{Q}_{2})$ defined as follows. Let $w=w_{\alpha_{i_{1}}}\ldots w_{\alpha_{i_{k}}}$, be a minimal expression for $w$ in terms of the generators $w_{\alpha_{1}},\ldots,w_{\alpha_{r}}$. Then we define $\textbf{s}((w,\lambda)) = w_{\alpha_{i_{1}}}(1)\ldots w_{\alpha_{i_{k}}}(1)2^{\lambda}$. This assignment can be shown to be independent of the minimal expression chosen for $w$ as follows. The Steinberg relations imply that the elements $w_{\alpha_{j}}(1)$ also satisfy the braid relations. Thus independence can be seen by adapting the proof of the Theorem in Section 29.4 in Humphreys \cite{H75}. 

Similarly, there is an (isomorphic) affine Weyl group $\widehat{W}_{\mathrm{aff}}$ associated to the root system $(\frac{1}{2}\Phi,\mathcal{E})$, which can be decomposed as $\widehat{W}_{\mathrm{aff}}\cong  W\ltimes 2Y$. This group has a set of Coxeter generators given by $w_{\alpha_{1}}\ldots, w_{\alpha_{r}}$, and $w_{\alpha_{0}}\stackrel{\mathrm{def}}{=}w_{\alpha_{0},-2}$.

The group $\widehat{W}_{\mathrm{aff}}$ can be extended by the finite abelian group $\widetilde{Y}/2Y$ to give the extended affine Weyl group $\widetilde{W}_{\mathrm{aff}}\stackrel{\mathrm{def}}{=}W\ltimes \widetilde{Y} \cong (\widetilde{Y}/2Y)\ltimes \widehat{W}_{\mathrm{aff}}$. Let $\mathfrak{D}_{0}=\{v\in\mathcal{E}|0<\langle\alpha,v\rangle<1\text{ for all }\alpha\in \frac{1}{2}\Phi^{+}\}$ and let $\Omega=\mathrm{Stab}(\mathfrak{D}_{0})$.

There is a length function $\ell: \widetilde{W}_{\mathrm{aff}}\cong \widetilde{Y}\rtimes W \rightarrow \mathbb{Z}_{\geq0}$, which can be computed using the formula of Proposition 1.23 in Iwahori-Matsumoto \cite{IM65}. Namely, 
\begin{equation}\label{length}
\ell(\lambda s)=\sum_{\alpha\in\frac{1}{2}\Phi^{+}\cap s\frac{1}{2}\Phi^{+}}|\langle\alpha,\lambda\rangle|+\sum_{\alpha\in\frac{1}{2}\Phi^{+}\cap s\frac{1}{2}\Phi^{-}}|\langle\alpha,\lambda\rangle+1|,
\end{equation}
where $\lambda\in \widetilde{Y}$ and $s\in W$.

We view $\widetilde{W}_{\mathrm{aff}}$ as a subgroup of $W_{\mathrm{aff}}$ via the natural inclusion $W\ltimes \widetilde{Y}\hookrightarrow W\ltimes Y$. 

For an element $(s,\lambda)\in W_{\text{aff}}$, we will  sometimes abuse notation and let $s2^\lambda $ denote a representative of $(s,\lambda)$ in either $N$ or $\widetilde{N}$. 

\subsection{Induction and Restriction}

If $J$ is a locally compact Hausdorff topological group, let $\delta_{J}$ be a modular character of $J$.

If $(\sigma, \mathscr{V})$ is a smooth $\widetilde{B}$-representation, then the  normalized induction functor $i_{\widetilde{G},\widetilde{T}}=\mathrm{Ind}_{\widetilde{B}}^{\widetilde{G}}$ takes the $\widetilde{B}$-representation $\sigma$ to the $\widetilde{G}$-representation 

\begin{equation*}
i_{\widetilde{G},\widetilde{T}}(\sigma)=\{f:\widetilde{G}\rightarrow \mathscr{V}|\,f\text{ is smooth, and } f(bg)=\delta_{\widetilde{B}}(b)^{1/2}\sigma(b)f(g)\text{ for all }b\in \widetilde{B}.\},
\end{equation*}
where $\widetilde{G}$ acts by right translation.

Suppose $(\pi, \mathscr{V})$ is a smooth $\widetilde{G}$-representation. Let $\mathscr{V}(U)=\mathrm{span}\{\pi(u)v-v|u\in U^{*},\,v\in \mathscr{V}\}$. The normalized (Jacquet) restriction functor $r_{\widetilde{T},\widetilde{G}}$ takes a $\widetilde{G}$-representation $\pi$ to the $\widetilde{T}$-representation $\mathscr{V}_{U}=\mathscr{V}/\mathscr{V}(U)$, where the $\widetilde{T}$ action is defined by
\begin{equation*}
r_{\widetilde{T},\widetilde{G}}(t)(v+\mathscr{V}(U))=\delta_{\widetilde{B}}^{-1/2}(t)\pi(t)v+\mathscr{V}(U).
\end{equation*}

\subsection{The Covering Torus}\label{CoveringTorus}

In this subsection we recall some facts from Loke-Savin \cite{LS10a} about the structure of $\widetilde{T}$ and the classification of its genuine irreducible representations. Let $T^{\diamond}\stackrel{\text{def}}{=}T(\mathbb{Z}/4\mathbb{Z})\cong T(\mathbb{Z})\cong Y\otimes\mu_{2}$ and let $\widetilde{T}^{1}(\mathbb{Q}_{2})\stackrel{\mathrm{def}}{=}T^{*}_{1}\Upsilon(Y)\mu_{2}$. Regarding the structure of $\widetilde{T}$, Loke-Savin (page 4908) show that $\widetilde{T}\cong (\widetilde{T}^{\diamond}\times \widetilde{T}^{1}(\mathbb{Q}_{2}))/\Delta(\mu_{2})$, where $\Delta$ embeds $\mu_{2}$ along the diagonal. Furthermore, they relate $\widetilde{T}^{1}(\mathbb{Q}_{2})$ to a tame covering torus. For a precise statement see Proposition 4.5 \cite{LS10b}.

One consequence of this decomposition is that every genuine representation of $\widetilde{T}$ is the tensor product of a genuine representation of $\widetilde{T}^{\diamond}$ and a genuine representation of $\widetilde{T}^{1}(\mathbb{Q}_{2})$. The group $\widetilde{T}^{\diamond}$ is a finite two-step nilpotent group and an irreducible genuine $\widetilde{T}^{\diamond}$-representation is determined by its central character (Loke-Savin \cite{LS10b}, page 4907). 

We can also say something about the representations of $\widetilde{T}^{1}(\mathbb{Q}_{2})$. Let $T_{1,8}^{*}$ be the subgroup of $\widetilde{T}$ generated by the elements of the form $\tilde{h}_{\alpha}(u)$, where $u\in1+8\mathbb{Z}_{2}$. An irreducible genuine representations of the group $\widetilde{T}^{1}(\mathbb{Q}_{2})/T_{1,8}^{*}$ is determined by its central character (Loke-Savin \cite{LS10b}, Proposition 4.5 and Proposition 4.3). Thus any irreducible genuine $\widetilde{T}$-representation in which $T_{1,8}^{*}$ acts trivially is determined by the action of $Z(\widetilde{T}^{\diamond})$ and $Z(\widetilde{T}^{1}(\mathbb{Q}_{2}))$. 

For the remainder of this paper, we fix a Weyl group invariant genuine irreducible $\widetilde{T}^{\diamond}$-representation $(\tau^{\diamond},E)$. For existence see Lemma 4.11 \cite{ABPTV07}, where our $\widetilde{T}^{\diamond}$ is an example of the group $M$. In fact each genuine irreducible representation of $\widetilde{T}^{\diamond}$ is Weyl group invariant.

On page 4910 \cite{LS10b}, Loke-Savin introduce a genuine character $\gamma_{2}:Z(\widetilde{T}^{1}(\mathbb{Q}_{2}))\rightarrow \mu_{2}$ that is the identity on $Z(\widetilde{T}^{1}(\mathbb{Q}_{2}))\cap T^{*}_{1}$ and $\Upsilon(\widetilde{Y})$. Let $V(\gamma_{2})$ be the representation of $\widetilde{T}^{1}(\mathbb{Q}_{2})$ that is associated with $\gamma_{2}$. Now for any unramified $\chi:T\rightarrow\mathbb{C}^{\times}$, we define the genuine $\widetilde{T}$-representation $(\sigma_{\chi}, i(\chi))$, where $i(\chi)\stackrel{\text{def}}{=}(\tau^{\diamond}\otimes V(\gamma_{2}))\otimes\chi$. We will also use this notation for the inflation of $i(\chi)$ to $\widetilde{B}$. An unramified principal series of $\widetilde{G}$ is a representation of the form $i_{\widetilde{G},\widetilde{T}}(\sigma_{\chi})$.

One important property of $i(\chi)$ is that for any $w\in W$ we have $i(\chi)^{w}\cong i(\chi^{w})$. This follows from the Weyl group invariance of $\tau^{\diamond}$ and $V(\gamma_{2})$ (Loke-Savin \cite{LS10a}, Corollary 5.2).


\section{Splitting}\label{splitting}

This section contains two important results, Theorem \ref{Gamma14Splitting} and Proposition \ref{SplitWeylConjugate}. Theorem \ref{Gamma14Splitting} states that there is a group homomorphism $S:\Gamma_1(4)\rightarrow \widetilde{G}$ such that $\mathrm{pr}\circ S=id_{\Gamma_{1}(4)}$. (i.e., S splits sequence (\ref{CentralExt}).) This result is necessary to define the Hecke-algebra $\mathcal{H}$. Proposition \ref{SplitWeylConjugate} states that $S$ satisfies an important technical property used to construct a basis for $\mathcal{H}$ (Proposition \ref{BasisWellDefined}).


We begin with a few preliminaries.
\begin{lemma}\label{StGroupLemma} 
Let $\alpha\in\Phi^{+}$ and $u,t\in\mathbb{Q}_{2}^{\times}$ such that $1+tu\neq0$. Then in St$(\Phi,\mathbb{Q}_2)$

\begin{equation}
x^{\prime}_{\alpha}(t)x^{\prime}_{-\alpha}(u)=\{1+tu,\frac{t}{1+tu}\}^{-1}x^{\prime}_{-\alpha}(\frac{u}{1+tu})h^{\prime}_{\alpha}(1+tu)x^{\prime}_{\alpha}(\frac{t}{1+tu}).
\end{equation}
\end{lemma}

\textbf{Proof:} This follows from the Steinberg relations. Alternatively, an equivalent identity is a consequence Proposition 2.7 b) in Stein \cite{S73}.\EndProof

\begin{corollary}\label{ImportantFactorization}Let $\lambda\in\widetilde{Y}$. Let $u,t\in\mathbb{Q}_{2}$ such that $\mathrm{val}_{2}(t)\geq \langle\lambda,\alpha\rangle$ and $\mathrm{val}_{2}(u)\geq \langle\lambda,-\alpha\rangle+2$. Then the following identity holds in $\widetilde{G}$:

\begin{equation*}
\tilde{x}_{\alpha}(t)\tilde{x}_{-\alpha}(-u)=\tilde{x}_{-\alpha}(\frac{u}{1+tu})\tilde{h}_{\alpha}(1+tu)\tilde{x}_{\alpha}(\frac{t}{1+tu}).
\end{equation*}
\end{corollary}
\textbf{Proof:} Note that $1+tu\in1+4\mathbb{Z}_{2}$, since $\mathrm{val}_{2}(t)\geq \langle\lambda,\alpha\rangle$ and $\mathrm{val}_{2}(u)\geq \langle\lambda,-\alpha\rangle+2$. Thus, $(1+tu,\frac{t}{1+tu})_{2}=(1+tu,t)_{2}$. Let $t=2^{\langle\lambda,\alpha\rangle}t^{\prime}$, where $t^{\prime}\in\mathbb{Z}_{2}$. Since $\lambda\in\widetilde{Y}$ we have $(1+tu,t)_{2}=(1+tu,t^\prime)_{2}$. Now $(1+tu,t^\prime)_{2}=1$, because $2|t^\prime$ implies $1+tu\equiv 1$ (mod $8$). \EndProof

Now we can prove that $\Gamma_{1}(4)$ splits sequence (\ref{CentralExt}). Let $\Gamma_{1}(4)^{*}$ be the subgroup of $\widetilde{G}$ generated by the elements $\tilde{x}_{\alpha}(t)$, $\tilde{x}_{-\alpha}(4u)$, $\tilde{h}_{\alpha}(v)$ for all $\alpha\in \Phi^{+}$ $t,u\in\mathbb{Z}_{2}$ and $v\in 1+4\mathbb{Z}_{2}$.

\begin{theorem}\label{Gamma14Splitting} The group homomorphism $\mathrm{pr}:\Gamma_{1}(4)^{*}\rightarrow \Gamma_{1}(4)$ is an isomorphism. Moreover, its inverse $S:\Gamma_{1}(4)\rightarrow \Gamma_{1}(4)^{*}$ splits sequence (\ref{CentralExt}).
\end{theorem}

\textbf{Proof:} The group $\Gamma_{1}(4)$ is generated by the elements $x_{\alpha}(t)$, $x_{-\alpha}(4u)$, and $h_{\alpha}(v)$ for all $\alpha\in \Phi^{+}$ $t,u\in\mathbb{Z}_{2}$ and $v\in 1+4\mathbb{Z}_{2}$. A complete set of relations for $\Gamma_{1}(4)$ is given by the Steinberg relations and the identity of Corollary \ref{ImportantFactorization}. The Steinberg relations and Corollary \ref{ImportantFactorization} also form a complete set of relations for the group $\Gamma_{1}(4)^{*}$. Since the projection map sends the generators of $\Gamma_{1}(4)^{*}$ to the generators of $\Gamma_{1}(4)$, this map is an isomorphism. The inverse map $S$ is a splitting by definition.\EndProof 

\begin{proposition}\label{SplitWeylConjugate}
Let $x\in \widetilde{N}^{\prime}$. Then $\widetilde{\Gamma}_{0}(4)\cap x\Gamma_{1}(4)^{*}x^{-1}\subseteq \Gamma_{1}(4)^{*}$.
\end{proposition}

\textbf{Proof:} Let $\beta_{1},\ldots, \beta_{\ell}$ be an enumeration of the positive roots and consider the element 
\begin{equation}
\gamma=\tilde{x}_{-\beta_{\ell}}(u_{\ell})\ldots \tilde{x}_{-\beta_{1}}(u_{1})\tilde{x}_{\beta_{1}}(t_{1})\ldots \tilde{x}_{\beta_{\ell}}(t_{\ell})h,
\end{equation}
where $t_{i},u_{i}\in \mathbb{Q}_{2}$, and $h\in \widetilde{T}$. Then $\gamma\in \Gamma_{1}(4)^{*}$ if and only if $t_{i}\in \mathbb{Z}_{2}$, $u_{i}\in 4\mathbb{Z}_{2}$, and $h\in T_{1}$.
Furthermore, this factorization is unique.  The analogous statement holds if we permute the factors in any order or if we replace $\Gamma_{1}(4)^{*}$ by $\widetilde{\Gamma}_{0}(4)$. We will use these facts to prove the proposition.

Let $x=w2^{\lambda}$, where $w\in\mathcal{W}$ and $\lambda\in \widetilde{Y}$. Suppose that $x\gamma x^{-1}\in \widetilde{\Gamma}_{0}(4)$. Thus 
\begin{multline}\label{ConjGamma}
x\gamma x^{-1}=\tilde{x}_{-w\cdot\beta_{\ell}}(\pm2^{\langle\lambda,-\beta_{\ell}\rangle}u_{\ell})\ldots \tilde{x}_{-w\cdot\beta_{1}}(\pm2^{\langle\lambda,-\beta_{1}\rangle}u_{1})\\
\times \tilde{x}_{w\cdot\beta_{1}}(\pm2^{\langle\lambda,\beta_{1}\rangle}t_{1})\ldots \tilde{x}_{w\cdot\beta_{\ell}}(\pm2^{\langle\lambda,\beta_{\ell}\rangle}t_{\ell})(xhx^{-1}).
\end{multline}
Note that $xhx^{-1}\in T^{*}_{1}$. 

By the unique factorization we see that the argument of $\tilde{x}_{\pm w\cdot\beta_{i}}$ is an element of $\mathbb{Z}_{2}$, if $\pm w\cdot\beta_{i}$ is positive, and $4\mathbb{Z}_{2}$, if $\pm w\cdot\beta_{i}$ is negative. Thus each factor in expression (\ref{ConjGamma}) is an element of $\Gamma_{1}(4)^{*}$ thus $x\gamma x^{-1}\in \Gamma_{1}(4)^{*}$. \EndProof

We close this section by identifying an obstruction which prevents the subgroup $\bold{G}(\mathbb{Z}_{2})$ from spliting the sequence (\ref{CentralExt}). Suppose that $S^{\prime}:\bold{G}(\mathbb{Z}_{2})\rightarrow \widetilde{G}$ splits sequence (\ref{CentralExt}). Then it follows that for any $\alpha\in\Phi$, we have $S^{\prime}(h_{\alpha}(-1))=\pm\tilde{h}_{\alpha}(-1)$. The element $h_{\alpha}(-1)\in \bold{G}(\mathbb{Z}_{2})$ has order $2$; the element $\pm\tilde{h}_{\alpha}(-1)\in\widetilde{G}$ has order $4$, since $(-1,-1)_{\mathbb{Q}_{2}}=-1$. Thus $S^{\prime}$ cannot exist. More generally, any subgroup of $\widetilde{G}$ which contains $h_{\alpha}(-1)$ cannot split the sequence (\ref{CentralExt}). Thus, the Iwahori subgroup $\Gamma_{0}(2)$ also does not split the sequence (\ref{CentralExt}). 

This obstruction does not appear in the tame case (i.e., GCD$(n,p)=1$), because the tame Hilbert Symbol of a local field $F$ is trivial on $\mathcal{O}_{F}^{\times}\times \mathcal{O}_{F}^{\times}$.


\section{$\Gamma_0(4)$ double cosets}\label{doublecosets}

In this section we compute representatives of the double coset space $\Gamma_0(4)\backslash G/ \Gamma_0(4)$ for the purpose of studying the $\widetilde{\Gamma}_{0}(4)$-equivariant Hecke Algebra $\mathcal{H}$. We begin with some notation. For $A\subseteq \Phi^{-}$ let $x_{A}(2)=\prod_{\alpha\in A}x_{\alpha}(2)$, where the product is taken with respect to some ordering of the elements of $A$. (Lemma \ref{singleCosets} shows that the choice of an order is immaterial to the study of $\Gamma_{0}(4)$ double cosets.) Now we can describe representatives of $\Gamma_0(4)\backslash G/ \Gamma_0(4)$.

\begin{proposition}\label{CosetUniRed} For every $g\in G$, there exists sets $A,B\subseteq \Phi^{-}$, $w\in W$, and $\lambda\in Y$ such that  

\begin{enumerate}
\item $\Gamma_0(4)g\Gamma_0(4)=\Gamma_0(4)x_{A}(2)w2^\lambda x_{B}(2)\Gamma_0(4)$;
\item $A\cap B=\emptyset$;\label{UniRed0}
\item $\langle\lambda,w^{-1}\alpha\rangle>0$ for all $\alpha\in A$; \label{UniRed1}
\item$\langle\lambda,\beta\rangle\leq 0$ for all $\beta\in B$;\label{UniRed2}
\item for all $\beta\in B$, $\langle\lambda,\beta\rangle=0$ implies that $w\beta\in \Phi^{-}$.\label{UniRed3}
\end{enumerate}
\end{proposition}

We will prove this proposition in a few steps. First observe that the group $\Gamma_0(4)$ is contained in the Iwahori subgroup $\Gamma_0(2)$. The Iwahori-Bruhat decomposition states that $\Gamma_0(2)\backslash G/ \Gamma_0(2)$ is in bijection with the affine Weyl group $W_{\text{aff}}\cong W\ltimes Y$. Thus to find a set of representatives of the doubles cosets of $G$ with respect to $\Gamma_0(4)$ we should determine representatives of $\Gamma_0(4)\backslash \Gamma_0(2)$ and $\Gamma_0(2)/\Gamma_0(4)$. 

\begin{lemma} \label{singleCosets}
The set $\{x_{A}(2)|A\subseteq \Phi^{-}\}$ is a complete set of distinct representatives of $\Gamma_0(4)\backslash \Gamma_0(2)$. The same holds for $\Gamma_0(2)/\Gamma_0(4)$. Furthermore, this holds for any permutation of the factors of $x_{A}(2)=\prod_{\alpha\in A}x_{\alpha}(2)$.
\end{lemma}

\textbf{Proof:} First we will show that $\Gamma_0(2)/\Gamma_0(4)$ is in bijection with $(\Gamma_0(2)\cap U_{-})/(\Gamma_0(4)\cap U_{-})$. The group $\Gamma_{0}(2)$ possesses an Iwahori factorization, thus multiplication defines a bijection 

\begin{equation*}
(\Gamma_{0}(2)\cap U_{-})\times(\Gamma_{0}(2)\cap B)\cong \Gamma_{0}(2).
\end{equation*}
Since reduction mod $2$ takes elements of $\textbf{G}(\mathbb{Z}_2)\cap B$ to $\textbf{B}(\mathbb{Z}/2\mathbb{Z})$, it follows that $\textbf{G}(\mathbb{Z}_2)\cap B=\Gamma_{0}(2)\cap B$. Thus

\begin{equation*}
(\Gamma_{0}(2)\cap U_{-})\times(\textbf{G}(\mathbb{Z}_2)\cap B)\cong \Gamma_{0}(2).
\end{equation*}
Similarly, 
\begin{equation*}
(\Gamma_{0}(4)\cap U_{-})\times(\textbf{G}(\mathbb{Z}_2)\cap B)\cong \Gamma_{0}(4).
\end{equation*}
Thus there is a bijection
\begin{equation*}
\Gamma_0(2)/\Gamma_0(4)\leftrightarrow(\Gamma_0(2)\cap U_{-})/(\Gamma_0(4)\cap U_{-}).
\end{equation*}

Since $U_{-}$ is a smooth group scheme over $\mathbb{Z}$ it is a smooth group scheme over $\mathbb{Z}_2$. This implies that $(\Gamma_0(2)\cap U_{-})/(\Gamma_0(4)\cap U_{-})$, which is an abelian group, can be identified with $\mathfrak{u}_{-}(\mathbb{F}_2)$, the Lie algebra of $\textbf{U}_{-}$ over the field with two elements, as abelian groups. In particular, this map sends the Chevalley generators $X_\alpha$ of $\mathfrak{u}_{-}(\mathbb{F}_2)$ to $x_{\alpha}(2)$. Since $\mathfrak{u}_{-}(\mathbb{F}_2)$ is abelian we see that the order of the elements in $A$ does not change the coset. The analogous argument proves the result for $\Gamma_0(4)\backslash \Gamma_0(2)$.\EndProof

Lemma \ref{singleCosets} implies that a complete set of representatives of $\Gamma_0(4)\backslash G/ \Gamma_0(4)$ is contained among the elements of the set 
\[\{x_{A}(2)w2^\lambda x_{B}(2)| w\in W,\,A,B\subseteq \Phi^{-},\,\lambda\in Y\}.\]
The $x_{A}(2)$ will be referred to as the unipotent elements on the left and the $x_{B}(2)$ will be referred to as the unipotent elements on the right.

Our next task is to eliminate redundant representatives. However, we will stop short of finding a complete set of distinct representatives. 

\begin{lemma}\label{removeuni} \label{LRrelated} Let $A,B\subseteq \Phi^{-}$, $\alpha\in A$, $\beta\in B$, $w\in W$, and $\lambda\in Y$. Let $A^\prime=A-\{\alpha\}$ and $B^\prime=B-\{\beta\}$.
\begin{enumerate}
\item\label{RemoveLeft} If $\langle\lambda,-w^{-1}\alpha\rangle\geq 1$, or if $\langle\lambda,w^{-1}\alpha\rangle=0$ and $w^{-1}\cdot\alpha\in\Phi^{+}$, then 
\begin{equation}
\Gamma_0(4)x_{A}(2)w2^\lambda x_{B}(2)\Gamma_0(4)=\Gamma_0(4)x_{A^\prime}w2^\lambda x_{B}(2)\Gamma_0(4).
\end{equation}

\item\label{RemoveRight} If $\langle\lambda,\beta\rangle\geq 1$, or if $\langle\lambda,\beta\rangle=0$ and $w\cdot\beta\in\Phi^{+}$, then 
\begin{equation}
\Gamma_0(4)x_{A}(2)w2^\lambda x_{B}(2)\Gamma_0(4)
=\Gamma_0(4)x_{A}(2)w2^\lambda x_{B^\prime}(2)\Gamma_0(4).
\end{equation}

\item\label{LtoR} If $\langle\lambda,w^{-1}\alpha\rangle=0$ and $w^{-1}\cdot\alpha\in\Phi^{-}$, then 
\begin{equation}
\Gamma_0(4)x_{A}(2)w2^\lambda x_{B}(2)\Gamma_0(4)\\
=\Gamma_0(4)x_{A^{\prime}}(2)w2^\lambda x_{B\cup\{w^{-1}\cdot\alpha\}}(2)\Gamma_0(4).
\end{equation}

\item\label{RtoL} If $\langle\lambda,\beta\rangle=0$ and $w\cdot\beta\in\Phi^{-}$, then 
\begin{equation}
\Gamma_0(4)x_{A}(2)w2^\lambda x_{B}(2)\Gamma_0(4)\\
=\Gamma_0(4)x_{A\cup\{w\cdot\beta\}}(2)w2^\lambda x_{B^{\prime}}(2)\Gamma_0(4).
\end{equation}
\end{enumerate}
\end{lemma}

\textbf{Proof:}  We will prove statement (\ref{RemoveLeft}); the proofs of the remaining statements are identical. 

We have $\Gamma_{0}(4)x_{A}(2)=\Gamma_{0}(4)x_{A^\prime}(2)x_{\alpha}(2)$, and $x_{\alpha}(2)w2^\lambda=w2^\lambda x_{w^{-1}\alpha}(\pm2^{1+\langle\lambda,-w^{-1}\alpha\rangle})$. Thus 
\begin{align*}
\Gamma_0(4)x_{A}(2)w2^\lambda x_{B}(2)\Gamma_0(4)=&\Gamma_0(4)x_{A^\prime}(2)x_{\alpha}(2)w2^\lambda x_{B}(2)\Gamma_0(4)\\
=&\Gamma_0(4)x_{A^\prime}(2)w2^\lambda x_{w^{-1}\alpha}(\pm2^{1+\langle\lambda,-w^{-1}\alpha\rangle})x_{B}(2)\Gamma_0(4).
\end{align*}

Since $\langle\lambda,-w^{-1}\alpha\rangle\geq 1$, we have $x_{w^{-1}\alpha}(\pm2^{1+\langle\lambda,-w^{-1}\alpha\rangle})\in\Gamma(4)$. The subgroup $\Gamma(4)$ is normal in $G(\mathbb{Z})$ and contained in $\Gamma_{0}(4)$ so the element $x_{w^{-1}\alpha}(\pm2^{1+\langle\lambda,-w^{-1}\alpha\rangle})$ can be moved right and absorbed into $\Gamma_0(4)$.  \EndProof\\

\textbf{Proof of Proposition \ref{CosetUniRed}:} We proved that $\Gamma_{0}(4)g\Gamma_{0}(4)=\Gamma_0(4)x_{A^{\prime}}(2)w2^\lambda x_{B^{\prime}}(2)\Gamma_0(4)$ for some $A^{\prime},B^{\prime}\subseteq \Phi^{-}$, $w\in W$, and $\lambda\in Y$. Using Lemma \ref{removeuni} we can identify a preferred representative of $\Gamma_0(4)x_{A^{\prime}}(2)w2^\lambda x_{B^{\prime}}(2)\Gamma_0(4)$. For each $\alpha\in A^{\prime}$ we can check to see if $\alpha$ satisfies the hypotheses of item (\ref{RemoveLeft}) in Lemma \ref{removeuni}, in which case $\alpha$ can be removed without changing the double coset. Futhermore, if $\alpha$ satisfies item (\ref{LtoR}), then we can move $\alpha$ to the right-hand side where it becomes $w^{-1}\alpha$. Let $A$ be the set of elements in $A^{\prime}$ that do not satisfy the hypotheses of items (\ref{RemoveLeft}), and (\ref{LtoR}). Similarly, for each $\beta\in B^{\prime}$ we can check to see if $\beta$ satisfies the hypotheses of item (\ref{RemoveRight}) in Lemma \ref{removeuni}, in which case $\beta$ can be removed. Let $B$ be the the set of elements of $B^{\prime}$ that do not satisfy (\ref{RemoveRight}) and the elements $w^{-1}\alpha$, where $\alpha\in A^{\prime}$ satisfies $(\ref{LtoR})$. Proposition \ref{removeuni} above proves that $\Gamma_0(4)x_{A}(2)w2^\lambda x_{B}(2)\Gamma_0(4)=\Gamma_0(4)x_{A^\prime}(2)w2^\lambda x_{B^\prime}(2)\Gamma_0(4)$. By construction the pair $(A,B)$ satisfy the conditions of the proposition.\EndProof


\section{The $\Gamma_{0}(4)$ Hecke Algebra}\label{LinearHeckeAlg}

In this section we will study $\underline{\mathcal{H}}\stackrel{\mathrm{def}}{=}C_{c}^{\infty}(\Gamma_{0}(4)\backslash G/\Gamma_{0}(4))$, the algebra of $\Gamma_{0}(4)$-biinvariant compactly supported functions. The multiplication of $\underline{\mathcal{H}}$ is given by convolution $f_{1}*f_{2}(g)=\int_{G}f_{1}(h)f_{2}(h^{-1}g)dh$, where the Haar measure is normalized so that the measure of $\Gamma_{0}(4)$ is equal to $1$. 

To begin we introduce the following length function $\ell_{4}:G\rightarrow \mathbb{R}_{\geq 0}$ defined by
 \begin{equation}
 2^{\ell_{4}(g)} = [\Gamma_{0}(4):\Gamma_{0}(4)\cap g\Gamma_{0}(4)g^{-1}].
 \end{equation}
In fact, $\ell_{4}(\gamma_{1}g\gamma_{2})=\ell_{4}(g)$ for any $\gamma_{1},\gamma_{2}\in \Gamma_{0}(4)$. This can be seen as follows. The map $\Gamma_{0}(4)/\Gamma_{0}(4)\cap g\Gamma_{0}(4)g^{-1}\rightarrow \Gamma_{0}(4)g\Gamma_{0}(4)/\Gamma_{0}(4)$ defined by $\gamma\Gamma_{0}(4)\cap g\Gamma_{0}(4)g^{-1}\mapsto \gamma g\Gamma_{0}(4)$ is a bijection. Moreover, $ [\Gamma_{0}(4):\Gamma_{0}(4)\cap g\Gamma_{0}(4)g^{-1}]=|\Gamma_{0}(4)g\Gamma_{0}(4)/\Gamma_{0}(4)|$. 

Using the section $\textbf{s}:W_{\mathrm{aff}}\rightarrow G$ we can consider the function $\ell_{4}\circ\textbf{s}$. We will abuse notation and simply write $\ell_{4}(w)$ in place of $\ell_{4}\circ\textbf{s}(w)$. The main result of this section, Theorem \ref{EqualLengths}, states that $\ell_{4}$ and $\ell$ are proportional when restricted to $\widetilde{W}_{\mathrm{aff}}$. We begin with some preliminary results.

\begin{lemma} \label{DoubleCosetSpit}
Let $w_{1},w_{2}\in W_{\mathrm{aff}}$. Then $\Gamma_{0}(4)\textbf{s}(w_{1}w_{2})\Gamma_{0}(4)=\Gamma_{0}(4)\textbf{s}(w_{1})\textbf{s}(w_{2})\Gamma_{0}(4)$.
\end{lemma}
\textbf{Proof:} This follows because $\textbf{s}(w_{1}w_{2})^{-1}\textbf{s}(w_{1})\textbf{s}(w_{2})\in \textbf{T}(\mathbb{Z}_{2})\subset \Gamma_{0}(4)$.\EndProof

\begin{proposition} \label{ell2Finite}
For any $g\in G$, $\ell_{4}(g)<\infty$.
\end{proposition}
 
\textbf{Proof:} It suffices to show that the number of cosets $\delta\Gamma_{0}(4)\subset \Gamma_{0}(4)g\Gamma_{0}(4)$ is finite. By Proposition 3.1 in Iwahori-Matsumoto \cite{IM65} we know that the number of cosets $\delta\Gamma_{0}(2)\subset \Gamma_{0}(2)g\Gamma_{0}(2)$ is finite. Since $[\Gamma_{0}(2):\Gamma_{0}(4)]<\infty$ the result follows. \EndProof
 
Let $\mathrm{ind}:\underline{\mathcal{H}}\rightarrow \mathbb{C}$ be the algebra homomorphism defined by $f\mapsto \int_{G}f(h)dh$. For the characteristic function $\underline{e}_{g}=\mathds{1}_{\Gamma_{0}(4)g\Gamma_{0}(4)}$ we have $\mathrm{ind}(e_{g})=|\Gamma_{0}(4)g\Gamma_{0}(4)/\Gamma_{0}(4)|= 2^{\ell_{4}(g)}$.

\begin{proposition}\label{lengthmultLinear}
Let $g_{1},g_{2}\in G$. If $\ell_{4}(g_{1}g_{2})=\ell_{4}(g_{1})+\ell_{4}(g_{2})$, then $\underline{e}_{g_{1}}*\underline{e}_{g_{2}}=\underline{e}_{g_{1}g_{2}}$.
\end{proposition}

\textbf{Proof:} By definition 

\begin{equation}\label{ConvoSum}
\underline{e}_{g_{1}}*\underline{e}_{g_{2}}(g)=\int_{G}\underline{e}_{g_{1}}(h)\underline{e}_{g_{2}}(h^{-1}g)=\sum_{\delta\in G/\Gamma_{0}(4)}\underline{e}_{g_{1}}(\delta)\underline{e}_{g_{2}}(\delta^{-1}g).
\end{equation}
The right hand side of equation (\ref{ConvoSum}) is equal to the number of cosets $\delta\Gamma_{0}(4)$ satisfying 
\begin{equation}\label{Gamma04CosetConditions}
\delta\Gamma_{0}(4)\subset \Gamma_{0}(4)g_{1}\Gamma_{0}(4)\hspace{1cm} \text{ and } \hspace{1cm} (\delta\Gamma_{0}(4))^{-1}g\subseteq \Gamma_{0}(4)g_{2}\Gamma_{0}(4).
\end{equation}
First we consider the case where $g=g_{1}g_{2}$. In this case we can directly check that $g_{1}\Gamma_{0}(4)$ is one such coset. Thus it follows that $\underline{e}_{g_{1}}*\underline{e}_{g_{2}}=c\underline{e}_{g_{1}g_{2}}+f$, where $c\in\mathbb{Z}_{\geq 1}$ and $f\in \underline{\mathcal{H}}$ is a nonnegative function.

We can apply $\mathrm{ind}$ to get $2^{\ell_{4}(g_{1}g_{2})}=c2^{\ell_{4}(g_{1}g_{2})}+\mathrm{ind}(f)$, because $\ell_{4}(g_{1}g_{2})=\ell_{4}(g_{1})+\ell_{4}(g_{2})$. Since $\mathrm{ind}(f)\geq 0$ and $c\geq 1$ it follows that $f=0$ and $c=1$. Thus $\underline{e}_{g_{1}}*\underline{e}_{g_{2}}=\underline{e}_{g_{1}g_{2}}$. \EndProof

\textbf{Remark:} The above proof also shows that if $\delta\Gamma_{0}(4)$ satisfies the coset conditions of (\ref{Gamma04CosetConditions}), then $\Gamma_{0}(4)g\Gamma_{0}(4)=\Gamma_{0}(4)g_{1}g_{2}\Gamma_{0}(4)$ and $\delta\Gamma_{0}(4)=g_{1}\Gamma_{0}(4)$.

\begin{proposition}
Let $w=w_{\alpha_{i}}\in\widetilde{W}_{\mathrm{aff}}$, where $i=0,\ldots, r$. Then 
\begin{equation}
\mathrm{ind}(\underline{e}_{\textbf{s}(w)})=2^{2}.
\end{equation}
\end{proposition}
  
\textbf{Proof:} If $i\neq 0$, then the Iwahori factorization implies that 
\begin{equation}
\Gamma_{0}(4)w_{\alpha_{i}}(1)\Gamma_{0}(4)/\Gamma_{0}(4)\cong U_{\alpha_{i}}(\mathbb{Z}_{2})/U_{\alpha_{i}}(4\mathbb{Z}_{2})\cong \mathbb{Z}/4\mathbb{Z}.
\end{equation}
  
If $i=0$, then the Iwahori factorization implies that 
\begin{equation}\label{AffineDoubleCoset}
\Gamma_{0}(4)w_{\alpha_{0}}(4)\Gamma_{0}(4)/\Gamma_{0}(4)\cong U_{\alpha_{0}}(4\mathbb{Z}_{2})/U_{\alpha_{0}}(4^{2}\mathbb{Z}_{2})\cong \mathbb{Z}/4\mathbb{Z}.
\end{equation}
\EndProof

\begin{lemma}\label{StabLen}
Let $w_{1}\in \Omega$ and $w_{2}\in \widetilde{W}_{\mathrm{aff}}$. Then $\ell_{4}(w_{1}w_{2})=\ell_{4}(w_{2})$.
\end{lemma}

\textbf{Proof:} First we claim that $\textbf{s}(w_{1})\in N_{G}(\Gamma_{0}(4))$. By Proposition 1.10. in Iwahori-Matsumoto \cite{IM65}, $w_{1}\in\Omega$ implies that $\ell(w_{1})=0$. Using formula (\ref{length}), the Iwahori factorization, and the Steinberg relations one can directly check that $\textbf{s}(w_{1})\in N_{G}(\Gamma_{0}(4))$. 

Finally, by Lemma \ref{DoubleCosetSpit} and $\Gamma_{0}(4)\textbf{s}(w_{1})\textbf{s}(w_{2})\Gamma_{0}(4)=\textbf{s}(w_{1})\Gamma_{0}(4)\textbf{s}(w_{2})\Gamma_{0}(4)$, we see that $\ell_{4}(w_{1}w_{2})=\ell_{4}(w_{2})$.\EndProof

Now we can prove the main theorem of this section.
 
\begin{theorem}\label{EqualLengths}
Let $w \in \widetilde{W}_{\mathrm{aff}}$. Then 
\begin{equation}
\ell_{4}(w)=2\ell(w).
\end{equation}
\end{theorem}
 
\textbf{Proof:} Suppose that $w=w_{1}w_{2}$, where $w_{1}\in \Omega$ and $w_{2}\in \widehat{W}_{\mathrm{aff}}$. Then $\ell(w_{1}w_{2})=\ell(w_{2})$. By Lemma \ref{StabLen}, $\ell_{4}(w_{1}w_{2})=\ell_{4}(w_{2})$. Thus it suffices to show that $\ell_{4}(w)=2\ell(w)$ for $w\in \widehat{W}_{\mathrm{aff}}$.

Suppose that $w=w_{\alpha_{i_{1}}}\ldots w_{\alpha_{i_{k}}}$ is a minimal expression of $w$ with respect to the Coxeter generators $w_{\alpha_{0}}, \ldots w_{\alpha_{r}}$ of $\widehat{W}_{\mathrm{aff}}$. By Lemma \ref{DoubleCosetSpit} and Proposition \ref{lengthmultLinear}
\begin{equation}
2^{\ell_{4}(w)}=\mathrm{ind}(\underline{e}_{\textbf{s}(w)})=\mathrm{ind}(\underline{e}_{\textbf{s}(w_{\alpha_{i_{1}}})})\ldots\mathrm{ind}(\underline{e}_{\textbf{s}(w_{\alpha_{i_{k}}})})=2^{2k}=2^{2\ell(w)}.
\end{equation}
Thus $\ell_{4}(w)=2\ell(w)$ for $w\in \widehat{W}_{\mathrm{aff}}$.\EndProof 


\section{The $\widetilde{\Gamma}_{0}(4)$ Hecke Algebra} \label{2adicHecke}

Let $(\tau^{\diamond},E)$ be a finite dimensional irreducible genuine Weyl group-invariant representation of $\widetilde{T}^{\diamond}$. Since $\widetilde{\Gamma}_{0}(4)\cong \widetilde{T}^{\diamond}\ltimes S(\Gamma_{1}(0))$, $\tau^{\diamond}$ inflates to a representation of $\widetilde{\Gamma}_{0}(4)$ which we call $\tau$. Let $\tau^{\vee}$ be the contragredient of $\tau$. The $\tau$-spherical Hecke algebra of $\widetilde{G}$ is 
\[\mathcal{H}\stackrel{\mathrm{def}}{=}\mathcal{H}(\widetilde{G},\tau^{\vee})=\{f\in C^{\infty}_{c}(\widetilde{G},\text{End}(E))| f( k_1gk_2)=\tau(k_1)f(g)\tau(k_2),\,\text{for all }k_1,k_2\in \widetilde{\Gamma}_0(4)\}.\]
For $f_{1},f_{2}\in \mathcal{H}$, the multiplication is defined by $f_{1}*f_{2}(g)=\int_{\widetilde{G}}f_{1}(h)f_{2}(h^{-1}g)dh$, where the Haar measure on $\widetilde{G}$ is normalized so that $\widetilde{\Gamma}_{0}(4)$ has measure 1.

The main result of this section, Theorem \ref{HeckeIso}, describes a Bernstein presentation of $\mathcal{H}$. 

Now we outline our approach to Theorem \ref{HeckeIso}. First, we construct a $\mathbb{C}$-basis for $\mathcal{H}$, Proposition \ref{BasisWellDefined}. Second, we identify some multiplicative relations among these basis elements, propositions \ref{LengthMult} and \ref{QuadRel}. Third, we use propositions \ref{LengthMult} and \ref{QuadRel} to prove that $\mathcal{H}$ admits an Iwahori-Matsumoto presentation, Proposition \ref{AffineHeckePresentation}.  

Finally, using results of Lusztig \cite{L89}, the Iwahori-Matsumoto presentation implies that $\mathcal{H}$ satisfies the Bernstein relations, and Savin's trick (Lemma 7.6, \cite{S04}) shows that these relations imply all others. This results in Theorem \ref{HeckeIso}, a Bernstein presentation of $\mathcal{H}$.

Our first step is to construct a $\mathbb{C}$-basis for $\mathcal{H}$. After a few technical preliminaries, we show that the nontrivial action of $\mu_2$ prohibits certain double cosets from supporting any functions in $\mathcal{H}$, propositions \ref{NoUni} and \ref{ModCoCharRed}. After identifying these constraints on support, we can construct a basis for $\mathcal{H}$, Proposition \ref{BasisWellDefined}.
  
\begin{proposition}[Stein \cite{S73}, Corollary 2.9] \label{SteinCommutator}Let $\alpha\in\Phi$. Then $[\tilde{x}_{\alpha}(2),\tilde{x}_{-\alpha}(2)]=(-1)\gamma$, where $-1\in\mu_2$ and $\gamma\in S(\Gamma(4))$.
\end{proposition} 

\textbf{Proof:} This follows by applying Lemma \ref{StGroupLemma} and the Steinberg relations.\EndProof

\begin{lemma}\label{MetaMove}
Let $m\in\mathbb{Z}_{\geq 1}$, and let $\beta,\beta_j\in\Phi^-$. Then 
\[\tilde{x}_{\beta}(\pm2^m)\tilde{x}_{\beta_{1}}(2)\ldots \tilde{x}_{\beta_{\ell}}(2)S(\Gamma_1(4))=\tilde{x}_{\beta_{1}}(2)\ldots \tilde{x}_{\beta_{\ell}}(2)\tilde{x}_{\beta}(\pm2^m)\Gamma_1(4)^{*},\]
and
\[\Gamma_1(4)^{*}\tilde{x}_{\beta}(\pm2^m)\tilde{x}_{\beta_{1}}(2)\ldots \tilde{x}_{\beta_{\ell}}(2)=\Gamma_1(4)^{*}\tilde{x}_{\beta_{1}}(2)\ldots \tilde{x}_{\beta_{\ell}}(2)\tilde{x}_{\beta}(\pm2^m).\]

\end{lemma}

\textbf{Proof:} The result follows from induction on $\ell$ and the Steinberg relations. \EndProof

\begin{lemma}\label{MetaMove2}
Let $m\in\mathbb{Z}_{\geq 1}$, $\eta\in \Phi,\beta_{j}\in\Phi^-$, $w\in W$, and $\lambda\in X_*(T)$, such that $-w\cdot\eta\in \Phi^{+}$, $\langle\lambda,w^{-1}\cdot\beta_{j}\rangle>0$ for all $j$, and $\langle\lambda,-\eta\rangle\geq0$. Then 
\[\Gamma_1(4)^{*}\tilde{x}_{\beta_{1}}(2)\ldots \tilde{x}_{\beta_{k}}(2)\tilde{x}_{-w\cdot\eta}(\pm2^m)=\Gamma_1(4)^{*}\tilde{x}_{\beta_{1}}(2)\ldots \tilde{x}_{\beta_{k}}(2).\]

\end{lemma}

\textbf{Proof:} We prove this proposition using induction on $k$. When $k=0$ the result holds since $\tilde{x}_{-w\cdot\eta}(\pm2^m)\in \Gamma_{1}(4)^{*}$. Now suppose that the result holds for $k-1$.

Consider the identity $\tilde{x}_{\beta_{k}}(2)\tilde{x}_{-w\cdot\eta}(\pm2^m)=[\tilde{x}_{\beta_{k}}(2),\tilde{x}_{-w\cdot\eta}(\pm2^m)]\tilde{x}_{-w\cdot\eta}(\pm2^m)\tilde{x}_{\beta_{k}}(2)$. To compute the commutator we prove that $-\beta_{k}\neq -w\cdot\eta$. Suppose $-\beta_{k}= -w\cdot\eta$, then $0>\langle\lambda,-w^{-1}\cdot\beta_{k}\rangle=\langle\lambda,-\eta\rangle\geq 0$. This is a contradiction so $-\beta_{k}\neq -w\cdot\eta$.

Thus the commutator can be computed using the Steinberg relations. Since our root system is simply laced we see that if $-w\cdot\eta+\beta_{k}\notin\Phi$, then $[\tilde{x}_{\beta_{k}}(2),\tilde{x}_{-w\cdot\eta}(\pm2^m)]=1$; if $-w\cdot\eta+\beta_{k}\in\Phi$, then $[\tilde{x}_{\beta_{k}}(2),\tilde{x}_{-w\cdot\eta}(\pm2^m)]=\tilde{x}_{-w\cdot\eta+\alpha_k}(\pm2^{m+1})$. If $-w\cdot\eta+\alpha_k\in\Phi^{-}$ we can apply Lemma \ref{MetaMove} to remove $\tilde{x}_{-w\cdot\eta+\alpha_k}(\pm2^{m+1})$. Then we may apply induction to remove $\tilde{x}_{-w\cdot\eta}(\pm2^m)$.

It remains to consider the case were $-w\cdot\eta+\alpha_k=w\cdot(-\eta+w^{-1}\cdot\alpha_k)\in\Phi^{+}$. But $\langle\lambda,-\eta+w^{-1}\cdot\alpha_k\rangle\geq 0$ so we may apply the induction hypothesis to remove $\tilde{x}_{-w\cdot\eta+\alpha_k}(\pm2^{m+1})$ and use the induction hypothesis once more to remove $\tilde{x}_{-w\cdot\eta}(\pm2^m)$.\EndProof\\

Now we can prove our first constraint on the support of functions in $\mathcal{H}$.

\begin{proposition} \label{NoUni}Let $A,B\subseteq \Phi^{-}$ and let $w\in W$ and $\lambda\in Y$. Suppose $(A,B)$ satisfy the conditions (\ref{UniRed0}), (\ref{UniRed1}), (\ref{UniRed2}), (\ref{UniRed3}) of Proposition \ref{CosetUniRed} and either $A$ or $B$ is nonempty. If $f\in \mathcal{H}$, then 
\[f(\widetilde{\Gamma}_0(4)\tilde{x}_{A}(2)w2^\lambda \tilde{x}_{B}(2)\widetilde{\Gamma}_0(4))=0.\]

\end{proposition}

\textbf{Proof:} Let us suppose that $B=\{\beta_{1},\ldots,\beta_{\ell}\}$ is nonempty and 
\[\tilde{x}_{B}(2)=\tilde{x}_{\beta_{1}}(2)\ldots \tilde{x}_{\beta_{\ell}}(2).\]
Recall from Proposition \ref{CosetUniRed} that for all $j$
\begin{equation*}\langle\lambda,\beta_{j}\rangle\leq 0.
\end{equation*}
It will be convenient to single out $\beta_{\ell}$, so we will let $\eta=\beta_{\ell}$.

 The set $A$ may or may not be empty. If $A$ is not empty, Corollary \ref{CosetUniRed} implies for $\alpha\in A$
\[\langle\lambda,-w^{-1}\cdot\alpha\rangle<0.\]

Let $f\in\mathcal{H}$. We will prove that $f(\tilde{x}_{A}(2)w2^\lambda \tilde{x}_{B}(2))=-f(\tilde{x}_{A}(2)w2^\lambda \tilde{x}_{B}(2))$.

First note that $\tilde{x}_{-\eta}(2)\in S(\Gamma_1(4))$, so 
\begin{equation}f(\tilde{x}_{A}(2)w2^\lambda \tilde{x}_{\beta_{1}}(2)\ldots \tilde{x}_{\beta_{\ell}}(2))\\
=f(\tilde{x}_{A}(2)w2^\lambda \tilde{x}_{\beta_{1}}(2)\ldots \tilde{x}_{\beta_{\ell}}(2)\tilde{x}_{-\eta}(2)).\label{EQN1}
\end{equation}
By Proposition \ref{SteinCommutator} it follows that $\tilde{x}_{\beta_{\ell}}(2)\tilde{x}_{-\eta}(2)=\tilde{x}_{-\eta}(2)\tilde{x}_{\beta_{\ell}}(2)(-1)M$, where $M\in \Gamma_1(4)^{*}$. Therefore,
\begin{equation}
(\ref{EQN1})
=-f(\tilde{x}_{A}(2)w2^\lambda \tilde{x}_{\beta_{1}}(2)\ldots \tilde{x}_{\beta_{\ell-1}}(2)\tilde{x}_{-\eta}(2)\tilde{x}_{\beta_{\ell}}(2)).\label{EQN2}
\end{equation}
Now for any $j$ we have $[\tilde{x}_{\beta_{j}}(-2),\tilde{x}_{-\eta}(-2)]$ is equal to $1$ or $\tilde{x}_{\beta_{j}-\eta}(\pm2^{2})$. In either case $[\tilde{x}_{\beta_{j}}(-2),\tilde{x}_{-\eta}(-2)]\in S(\Gamma(4))$. These elements can be moved to the right since $S(\Gamma(4))$ is normalized by elements of the form $\tilde{x}_{\beta}(2)$, where $\beta\in\Phi$. By using induction on $\ell$ in conjunction with the Steinberg relations and Lemma \ref{MetaMove} it follows that 
\begin{equation}
(\ref{EQN2})\\
=-f(\tilde{x}_{A}(2)w2^\lambda \tilde{x}_{-\eta}(2)\tilde{x}_{B}(2)).\label{EQN3}
\end{equation}

We continue to push the element $\tilde{x}_{-\eta}(2)$ to the left using the identities $2^\lambda \tilde{x}_{-\eta}(2)=\tilde{x}_{-\eta}(2^{1+\langle\lambda,-\eta\rangle})2^\lambda$ and $w\tilde{x}_{-\eta}(2^{1+\langle\lambda,-\eta\rangle})=\tilde{x}_{-w\cdot\eta}(\pm2^{1+\langle\lambda,-\eta\rangle})w$. Therefore,
\begin{equation}
(\ref{EQN3})
=-f(\tilde{x}_{A}(2)\tilde{x}_{-w\cdot\eta}(\pm2^{1+\langle\lambda,-\eta\rangle})w2^\lambda \tilde{x}_{B}(2)).\label{EQN4}
\end{equation}

If $-w\cdot\eta\in\Phi^{-}$ and $\langle\lambda,\eta\rangle<0$, then we can apply Lemma \ref{MetaMove} to prove that 
\begin{equation}f(\tilde{x}_{A}(2)w2^\lambda \tilde{x}_{B}(2))
=-f(\tilde{x}_{A}(2)w 2^\lambda \tilde{x}_{B}(2)).
\end{equation}

If $-w\cdot\eta\in\Phi^{-}$ and $\langle\lambda,\eta\rangle=0$, then Proposition \ref{CosetUniRed} implies that $w\eta\in \Phi^{-}$, a contradiction. 

If $-w\cdot\eta\in\Phi^{+}$, then we can apply Lemma \ref{MetaMove2} to remove $\tilde{x}_{-w\cdot\eta}(\pm2^{1+\langle\lambda,-\eta\rangle})$. Thus we see that 
\begin{equation}f(\tilde{x}_{A}(2)wt(2^\lambda)\tilde{x}_{B}(2))\\
=-f(\tilde{x}_{A}(2)wt(2^\lambda)\tilde{x}_{B}(2)),
\end{equation}
as desired. 

One can apply a similar argument when $B$ is empty.\EndProof

\textbf{Remark:} The proof of Proposition \ref{NoUni} above does not utilize the full $\widetilde{\Gamma}_{0}(4)$ transformation law; it only requires $\Gamma_{1}(4)^{*}$-invariance.

To continue with the support calculations we will utilize the following lemma. Before we begin we introduce some notation. If $J$ is a subgroup of a group $K$ and $k\in K$, then let $^{k}J=kJk^{-1}$. If $\mathscr{V}$ is a representation of $J$ and $k\in K$, then we will write $^{k}\mathscr{V}$ for the representation of the group $kJk^{-1}$ acting on $\mathscr{V}$ via $(kjk^{-1})\cdot v=j\cdot v$.

\begin{lemma}[Bushnell-Kutzko \cite{BK93}, Proposition (4.1.1)]\label{BKsupp} Let $g\in \widetilde{G}$. The following are equivalent:\\
(i) Hom$_{\widetilde{\Gamma}_0(4)\cap ^g\widetilde{\Gamma}_{0}(4)}( ^g\mathscr{V},\mathscr{V})\neq0$;\\
(ii) there exists $\Phi\in \mathcal{H}$ with $\Phi(g)\neq 0$.

If the element $g$ satisfies these conditions, then we have a canonical vector space isomorphism between Hom$_{\widetilde{\Gamma}_0(4)\cap ^g\widetilde{\Gamma}_{0}(4)}( ^g\mathscr{V},\mathscr{V})$ and the space of functions $\Phi\in\mathcal{H}$ which vanish outside the double coset $\widetilde{\Gamma}_{0}(4)g\widetilde{\Gamma}_0(4)$.

\end{lemma}

Next we show that any double coset can support at most a one dimensional space of functions in $\mathcal{H}$.

\begin{proposition}\label{suppdim}
Let $w\in \mathcal{W}$, $\lambda\in Y$, and $g=w2^{\lambda}$, then 
\[\mathrm{dim}(\mathrm{Hom}_{\widetilde{\Gamma}_0(4)\cap ^g\widetilde{\Gamma}_{0}(4)}( ^gE,E))\leq 1.\] 
\end{proposition}

\textbf{Proof:} Since $\widetilde{T}^\diamond\subseteq \widetilde{\Gamma}_0(4)\cap ^g\widetilde{\Gamma}_{0}(4)$ we have Hom$_{\widetilde{\Gamma}_0(4)\cap\,^g\widetilde{\Gamma}_{0}(4)}( ^gE,E)\subseteq$ Hom$_{\widetilde{T}^\diamond}( ^gE,E)$. Recall that $(\tau^{\diamond},E)$ is Weyl group-invariant and irreducible, and note that the action of $2^\lambda$ on $\widetilde{T}^\diamond$ by conjugation is trivial. Therefore Schur's lemma implies that Hom$_{\widetilde{\Gamma}_0(4)\cap\,^g\widetilde{\Gamma}_{0}(4)}( ^gE,E)$ has dimension at most $1$.\EndProof

Now we can prove our second constraint on the support of functions in $\mathcal{H}$.

\begin{proposition} \label{ModCoCharRed}
Let $f\in\mathcal{H}$, $w\in \mathcal{W}$, and $\lambda\in Y-\widetilde{Y}$, then $f(w2^\lambda)=0$. \label{supp}
\end{proposition}

\textbf{Proof:} From Lemma \ref{BKsupp} it suffices to show that Hom$_{\widetilde{\Gamma}_0(4)\cap ^g\widetilde{\Gamma}_{0}(4)}( ^gE,E)=0$, where $g=w2^\lambda$. 
We begin with a few preliminary remarks. By assumption $\lambda\notin \widetilde{Y}$, so there exists $\alpha\in \Phi$ such that $\langle\lambda,\alpha\rangle$ is odd. Let $t\in \mathbb{Z}_{2}^{\times}$ such that $t\equiv 5$ ($8$). This implies that $(2,t)_{2}=-1$, $(-1,t)_2=1$, and $(t,t)_2=1$. Thus $h_{\alpha}(t)2^\lambda h_{\alpha}(t)^{-1}=-2^\lambda$, since $\langle\lambda,\alpha\rangle$ is odd; $wh_{\alpha}(t)w^{-1}=h_{w\alpha}(t)$, since $t\equiv 1$ (mod $4$).

Let $\psi\in \text{Hom}_{\widetilde{\Gamma}_0(4)\cap\,^g\widetilde{\Gamma}_{0}(4)}( ^gE,E)$. We will show that $\psi=-\psi$. For any $v\in E$, consider $\psi(w2^\lambda h_{\alpha}(t)(w2^\lambda)^{-1}v)$.  Since $\Gamma_{1}(4)^{*}$ acts trivially on $E$ and $t\equiv 1$ (mod $4$) it follow that that $h_{w\alpha}(t)$ acts trivially on $E$. Thus, $\psi(w2^\lambda h_{\alpha}(t)(w2^\lambda)^{-1}v)=-\psi(h_{w\alpha}(t)v)=-\psi(v)$. On the other hand, since $\psi$ is an intertwining operator and $t\equiv 1 $ (mod $4$) we have  $\psi(w2^\lambda h_{\alpha}(t)(w2^\lambda)^{-1}v)=h_{\alpha}(t)\psi(v)=\psi(v)$. Thus $\psi(v)=-\psi(v)$. Therefore, Hom$_{\widetilde{\Gamma}_0(4)\cap ^g\widetilde{\Gamma}_{0}(4)}( ^gE,E)=0$.\EndProof\\

Now we can construct a basis for $\mathcal{H}$. By Proposition 5.2 in Adams-Barbash-Paul-Trapa-Vogan \cite{ABPTV07}, the representation $\tau^{\diamond}$ can be extended from $\widetilde{T}^\diamond$ to $\mathcal{W}$. (Since the Hilbert symbol of $\mathbb{R}$ and $\mathbb{Q}_{2}$ agree on $\{\pm1\}\times\{\pm1\}$, the group $\mathcal{W}$ is a subgroup of the group $\widetilde{K}$ appearing in Proposition 5.2 of \cite{ABPTV07}.) We will call this extension $\tau_{\mathcal{W}}$. Recall that $\widetilde{N}^\prime\cong \mathcal{W}\ltimes \widetilde{Y}$. 
Thus $\tau_{\mathcal{W}}$ inflates to a representation of $\widetilde{N}^\prime$, which we call $\tau_{\widetilde{N}^{\prime}}$.

Proposition \ref{NoUni} and Proposition \ref{ModCoCharRed} together state that the support of a function in $\mathcal{H}$ is contained in the double cosets of the form $\widetilde{\Gamma}_{0}(4)x\widetilde{\Gamma}_{0}(4)$ where $x\in \widetilde{N}^\prime$. For each $w\in \widetilde{W}_{\mathrm{aff}}$ we define a function $e_{w}$. Let $x$ be any element of $\widetilde{N}^\prime$ that maps to $w$ under the natural map $\widetilde{N}^\prime\cong \mathcal{W}\ltimes \widetilde{Y}\rightarrow W\ltimes \widetilde{Y}\cong \widetilde{W}_{\mathrm{aff}}$. (Note that $\widetilde{T}^\diamond\subseteq \mathcal{W}$ is a normal subgroup, $\mathcal{W}/\widetilde{T}^{\diamond}\cong W$, and $\widetilde{T}^\diamond$ commutes with $\Upsilon(\widetilde{Y})$.) Define $e_{w}$ to be the unique function in $\mathcal{H}$ supported on $\widetilde{\Gamma}_{0}(4)x\widetilde{\Gamma}_{0}(4)$ such that 
\[e_{w}(\Gamma_1(4)^{*}x\Gamma_{1}(4)^{*})=\tau_{\widetilde{N}^{\prime}}(x).\] 
We must show that $e_{w}$ is well-defined and that the definition of $e_{w}$ is independent of our choice of $x$. Once we show that $e_{w}$ is well-defined, it is straightforward to show independence since every preimage of $w$ is of the form $xt$, where $t\in\widetilde{T}^{\diamond}$.

\begin{proposition}\label{BasisWellDefined}
The functions $e_{w}\in\mathcal{H}$, where $w\in \widetilde{W}_{\mathrm{aff}}$, are well-defined and form a $\mathbb{C}$-basis for $\mathcal{H}$.
\end{proposition}

\textbf{Proof:} (Well-defined) Suppose that $\gamma_{i}\in\widetilde\Gamma_{0}(4)$ and $x=\gamma_{1}x\gamma_{2}$. To prove that $e_{w}$ is well-defined we must prove that $\tau_{\widetilde{N}^{\prime}}(x)=\tau(\gamma_{1})\tau_{\widetilde{N}^{\prime}}(x)\tau(\gamma_{2})$. Let $\gamma_{1}=t_{1}u_{1}$ and $\gamma_{2}=u_{2}t_{2}$, where $t_{i}\in \widetilde{T}^{\diamond}$ and $u_{j}\in \Gamma_{1}(4)^{*}$. Thus, $\tau(\gamma_{1})\tau_{\widetilde{N}^{\prime}}(x)\tau(\gamma_{2})=\tau(t_{1})\tau_{\widetilde{N}^{\prime}}(x)\tau(t_{2})=\tau_{\widetilde{N}^{\prime}}(t_{1})\tau_{\widetilde{N}^{\prime}}(x)\tau_{\widetilde{N}^{\prime}}(t_{2})=\tau_{\widetilde{N}^{\prime}}(t_{1}xt_{2})$. Proposition \ref{SplitWeylConjugate} implies that $x=t_{1}xt_{2}$. Thus $\tau_{\widetilde{N}^{\prime}}(x)=\tau_{\widetilde{N}^{\prime}}(t_{1}xt_{2})$, as desired.

(Basis) Proposition \ref{suppdim} implies that $e_{w}$ generates the space of functions in $\mathcal{H}$ supported on the double coset $\widetilde{\Gamma}_{0}(4)x\widetilde{\Gamma}_{0}(4)$. Thus the functions $e_{w}$ form a basis of $\mathcal{H}$, as $w$ varies over elements of $\widetilde{W}_{\mathrm{aff}}\subseteq W_{\text{aff}}$.\EndProof


Our second step is to prove the multiplicative relations described in propositions \ref{LengthMult} and \ref{QuadRel}. Using these relations we prove that $\mathcal{H}$ admits an Iwahori-Matsumoto presentation in Proposition \ref{AffineHeckePresentation}.

\begin{proposition} \label{LengthMult}
If $w_1,w_2\in \widetilde{W}_{\mathrm{aff}}$ and $\ell(w_1w_2)=\ell(w_1)+\ell(w_2)$, then $e_{w_1}* e_{w_2}=e_{w_1w_2}$.
\end{proposition}

\textbf{Proof:} The proof is similar to the proof of Proposition 6.2 in Savin \cite{S04}. Let $x_{1}$, $x_{2}\in \widetilde{N}^{\prime}$ represent $w_{1}$ and $w_{2}$ respectively. For each $x\in\widetilde{N}^{\prime}$, to compute $e_{w_{1}}*e_{w_{2}}(x)$ we must find all cosets $\delta\widetilde{\Gamma}_{0}(4)$ such that 
\begin{equation}\label{nonlinearCosetCondition}
\delta\widetilde{\Gamma}_{0}(4)\subseteq \widetilde{\Gamma}_{0}(4)x_{1}\widetilde{\Gamma}_{0}(4)\hspace{1cm}\text{ and }\hspace{1cm} \widetilde{\Gamma}_{0}(4)\delta^{-1}x\subseteq \widetilde{\Gamma}_{0}(4)x_{2}\widetilde{\Gamma}_{0}(4).
\end{equation}
These containments only depend on $G$. So we may focus on finding all cosets $\mathrm{pr}(\delta)\Gamma_{0}(4)$ such that 
\begin{equation}\label{linearCosetCondition}
\mathrm{pr}(\delta)\Gamma_{0}(4)\subseteq \Gamma_{0}(4)\mathrm{pr}(x_{1})\Gamma_{0}(4)\hspace{.5cm}\text{ and }\hspace{.5cm} \Gamma_{0}(4)\mathrm{pr}(\delta)^{-1}\mathrm{pr}(x)\subseteq \Gamma_{0}(4)\mathrm{pr}(x_{2})\Gamma_{0}(4).
\end{equation}
By Theorem \ref{EqualLengths} and the remark after Proposition \ref{lengthmultLinear}, we know that $\Gamma_{0}(4)\mathrm{pr}(x)\Gamma_{0}(4)=\Gamma_{0}(4)\mathrm{pr}(x_{1}x_{2})\Gamma_{0}(4)$ and $\mathrm{pr}(\delta)\Gamma_{0}(4) = \mathrm{pr}(x_{1})\Gamma_{0}(4)$. Thus it is enough for us to compute $e_{w_{1}}*e_{w_{2}}(x_{1}x_{2})$. By definition
\begin{equation}\label{lengthmulteqn1}
e_{w_{1}}*e_{w_{2}}(x_{1}x_{2})=\int_{\widetilde{G}}e_{w_{1}}(h)e_{w_{2}}(h^{-1}x_{1}x_{2})dh=\sum_{\delta\in \widetilde{G}/\widetilde{\Gamma}_{0}(4)}e_{w_{1}}(\delta)e_{w_{2}}(\delta^{-1}x_{1}x_{2}).
\end{equation}
We have seen that the only coset that satisfies (\ref{nonlinearCosetCondition}) is $x_{1}\widetilde{\Gamma}_{0}(4)$. Thus
\begin{equation}\label{lengthmulteqn2}
(\ref{lengthmulteqn1}) = e_{w_{1}}(x_{1})e_{w_{2}}(x_{2})=e_{w_{1}w_{2}}(x_{1}x_{2}),
\end{equation}
as desired. The last equality follows because $\tau_{\widetilde{N}^{\prime}}$ is a representation of $\widetilde{N}^{\prime}$.\EndProof

For $\alpha\in\Delta\cup\{\alpha_{0}\}$ we will write $e_{\alpha}\stackrel{\text{def}}{=}e_{w_{\alpha}}$. Next we will prove that the elements $e_{\alpha}$ satisfy a quadratic relation. We begin with the following lemma.

\begin{lemma}\label{IrrationalWeyl} Let $\alpha\in \Phi$ and $k\in \mathbb{Z}$. Then 

\begin{equation}
\tau_{\widetilde{N}^{\prime}}(w_{\alpha}(2^{k}))+\tau_{\widetilde{N}^{\prime}}(w_{\alpha}(2^{k}))^{-1}=\epsilon \sqrt{2}I_{E},
\end{equation}
where $I_{E}$ is the identity map of $E$ and $\epsilon=\pm1$. Moreover, $\epsilon$ is independent of $\alpha$ and $k$.

\end{lemma}

\textbf{Proof:} In this proof we will write $w_{\alpha}$ for $\tilde{w}_{\alpha}(1)$. First, we show that $\tau_{\widetilde{N}^{\prime}}(w_{\alpha})+\tau_{\widetilde{N}^{\prime}}(w_{\alpha})^{-1}\in\mathrm{End}(E)$ is a scalar endomorphism. Specifically, for any $\beta\in \Delta$ we have 
\[\tau^{\diamond}(\tilde{h}_{\beta}(-1))[\tau_{\widetilde{N}^{\prime}}(w_{\alpha})+\tau_{\widetilde{N}^{\prime}}(w_{\alpha})^{-1}]\tau^{\diamond}(\tilde{h}_{\beta}(-1))^{-1}=\tau_{\widetilde{N}^{\prime}}(w_{\alpha})+\tau_{\widetilde{N}^{\prime}}(w_{\alpha})^{-1},\]
which follows from the Steinberg relations and because $\tau^{\diamond}(\tilde{h}_{\alpha}(-1))=\tau_{\widetilde{N}^{\prime}}(\tilde{h}_{\alpha}(-1))$. Since $(\tau^{\diamond},E)$ is irreducible, $\tau_{\widetilde{N}^{\prime}}(w_{\alpha})+\tau_{\widetilde{N}^{\prime}}(w_{\alpha})^{-1}$ is a scalar endomorphism, by Schur's Lemma.

So $\tau_{\widetilde{N}^{\prime}}(w_{\alpha})+\tau_{\widetilde{N}^{\prime}}(w_{\alpha})^{-1}=cI_{E}$ for some scalar $c$. We square this equation to get $\tau_{\widetilde{N}^{\prime}}(\tilde{h}_{\alpha}(-1))+\tau_{\widetilde{N}^{\prime}}(\tilde{h}_{\alpha}(-1))^{-1}+2I_{E}=c^2I_{E}$. However, $\tau_{\widetilde{N}^{\prime}}(\tilde{h}_{\alpha}(-1))+\tau_{\widetilde{N}^{\prime}}(\tilde{h}_{\alpha}(-1))^{-1}=\tau_{\widetilde{N}^{\prime}}(\tilde{h}_{\alpha}(-1))[I_{E}+\tau_{\widetilde{N}^{\prime}}(\tilde{h}_{\alpha}(-1))^{-2}]=\tau_{\widetilde{N}^{\prime}}(\tilde{h}_{\alpha}(-1))[I_{E}-I_{E}]=0$. Thus $c=\pm\sqrt{2}$.

Now consider $\tau_{\widetilde{N}^{\prime}}(\tilde{w}_{\alpha}(2^{k}))+\tau_{\widetilde{N}^{\prime}}(\tilde{w}_{\alpha}(2^{k}))^{-1}$. Since $\tilde{w}_{\alpha}(2^{k})=\tilde{h}_{\alpha}(2^{k})\tilde{w}_{\alpha}$ and $\tilde{w}_{\alpha}(-2^{k})=\tilde{h}_{\alpha}(2^{k})\tilde{w}_{\alpha}(-1)$ we see that
\begin{equation}
\tau_{\widetilde{N}^{\prime}}(\tilde{w}_{\alpha}(2^{k}))+\tau_{\widetilde{N}^{\prime}}(\tilde{w}_{\alpha}(2^{k}))^{-1} = \tau_{\widetilde{N}^{\prime}}(\tilde{h}_{\alpha}(2^{k}))[\epsilon\sqrt{2}I_{E}].
\end{equation}
But $\tau_{\widetilde{N}^{\prime}}(\tilde{h}_{\alpha}(2^{k}))=I_{E}$ so the result follows.

We have already seen that $\epsilon$ is independent of $k$. Now we will show that it is independent of $\alpha$. Since we are considering simply-laced root systems, the action of $W$ on $\Phi$ has a single orbit. Let $\alpha,\beta\in\Phi$ and let $w\in W$ such that $w(\alpha)=\beta$. Then by the Steinberg relations we have that 
\begin{equation}
\tau_{\widetilde{N}^{\prime}}(w)[\tau_{\widetilde{N}^{\prime}}(w_{\alpha})+\tau_{\widetilde{N}^{\prime}}(w_{\alpha})^{-1}]\tau_{\widetilde{N}^{\prime}}(w)^{-1}=\tau_{\widetilde{N}^{\prime}}(w_{w(\alpha)})+\tau_{\widetilde{N}^{\prime}}(w_{w(\alpha)})^{-1}=\tau_{\widetilde{N}^{\prime}}(w_{\beta})+\tau_{\widetilde{N}^{\prime}}(w_{\beta})^{-1}.
\end{equation}  
But since $\tau_{\widetilde{N}^{\prime}}(w_{\alpha})+\tau_{\widetilde{N}^{\prime}}(w_{\alpha})^{-1}$ is a scaler we also have 
\begin{equation}
\tau_{\widetilde{N}^{\prime}}(w)[\tau_{\widetilde{N}^{\prime}}(w_{\alpha})+\tau_{\widetilde{N}^{\prime}}(w_{\alpha})^{-1}]\tau_{\widetilde{N}^{\prime}}(w)^{-1}=\tau_{\widetilde{N}^{\prime}}(w_{\alpha})+\tau_{\widetilde{N}^{\prime}}(w_{\alpha})^{-1}.
\end{equation}
Thus $\epsilon$ is independent of  $\alpha$.
\EndProof 

\begin{proposition} \label{QuadRel}
Let $\alpha\in \Delta\cup\{\alpha_{0}\}$. Then 
\[e_{\alpha}^2=\epsilon\sqrt{2}e_{\alpha}+4\mathds{1},\]
where $\mathds{1}$ is the function in $\mathcal{H}$ supported on $\widetilde{\Gamma}_{0}(4)$ such that $\mathds{1}(\gamma)=\tau(\gamma)$.

\end{proposition}

\textbf{Proof:} We will address the case of $\alpha\in \Delta$ completely. When $\alpha=\alpha_{0}$ the computation is similar (because equation (\ref{AffineDoubleCoset}) holds) and will be omitted.

Suppose that $\alpha\in \Delta$. We claim that supp$(e_{\alpha}^2)\subseteq \widetilde{\Gamma}_{0}(4)\cup \widetilde{\Gamma}_{0}(4)w_\alpha\widetilde{\Gamma}_{0}(4)$. Using the Iwahori factorization and the Steinberg relations we can show that 
\begin{equation}
\widetilde{\Gamma}_{0}(4)w_{\alpha}(1)\widetilde{\Gamma}_{0}(4)w_{\alpha}(1)\widetilde{\Gamma}_{0}(4)=\cup_{t\in\mathbb{Z}/4\mathbb{Z}}\widetilde{\Gamma}_{0}(4)x_{-\alpha}(t)\widetilde{\Gamma}_{0}(4).
\end{equation}
For $t=0$, $\widetilde{\Gamma}_{0}(4)x_{-\alpha}(t)\widetilde{\Gamma}_{0}(4)=\widetilde{\Gamma}_{0}(4)$. For $t=\pm 1$, $\widetilde{\Gamma}_{0}(4)x_{\alpha}(t)\widetilde{\Gamma}_{0}(4)=\widetilde{\Gamma}_{0}(4)w_{\alpha}(1)\widetilde{\Gamma}_{0}(4)$. For $t=2$, $\widetilde{\Gamma}_{0}(4)x_{\alpha}(t)\widetilde{\Gamma}_{0}(4)=\widetilde{\Gamma}_{0}(4)x_{\alpha}(-2)w_{\alpha}(2)x_{\alpha}(-2)\widetilde{\Gamma}_{0}(4)$. Thus Proposition \ref{NoUni} implies supp$(e_{\alpha}^2)\subseteq \widetilde{\Gamma}_{0}(4)\cup \widetilde{\Gamma}_{0}(4)w_\alpha\widetilde{\Gamma}_{0}(4)$ and we can write $e_{\alpha}^2=ae_{\alpha}+b\mathds{1}$, where $a,b\in\mathbb{C}$.  To determine $a$ and $b$ it suffices to compute $e_{\alpha}^2(1)$ and $e_{\alpha}^2(w_{\alpha})$.

We begin with a preliminary computation.

\begin{align*}
e_{\alpha}^2(g)=&\int_{\widetilde{G}}e_{\alpha}(y)e_{\alpha}(y^{-1}g)dy\\
=&\sum_{y\in \widetilde{\Gamma}_{0}(4)w_{\alpha}\widetilde{\Gamma}_{0}(4)/\widetilde{\Gamma}_{0}(4)}e_{\alpha}(y)e_{\alpha}(y^{-1}g)\\
=&\sum_{\delta\in\widetilde{\Gamma}_{0}(4)/\widetilde{\Gamma}_{0}(4)\cap w_{\alpha}\widetilde{\Gamma}_{0}(4)w_{\alpha}^{-1}}e_{\alpha}(\delta w_{\alpha})e_{\alpha}(w_{\alpha}^{-1}\delta^{-1} g)\\
=&\sum_{u\in\mathbb{Z}/4\mathbb{Z}}e_{\alpha}(x_{\alpha}(u) w_{\alpha})e_{\alpha}(w_{\alpha}^{-1}x_{\alpha}(-u) g)\\
=&e_{\alpha}(w_{\alpha})\sum_{u\in\mathbb{Z}/4\mathbb{Z}}e_{\alpha}(w_{\alpha}^{-1}x_{\alpha}(u) g).
\end{align*}

It remains to compute $e_{\alpha}^2(1)$ and $e_{\alpha}^2(w_{\alpha})$. The computation of $e_{\alpha}^2(1)$ is straightforward and will be omitted. One finds that $e_{\alpha}^2(1)=4I=4\mathds{1}(1)$.

Finally we will compute $e_{\alpha}^2(w_{\alpha})$.

\begin{align*}
e_{\alpha}^2(w_{\alpha})=&e_{\alpha}(w_{\alpha})\sum_{u\in\mathbb{Z}/4\mathbb{Z}}e_{\alpha}(w_{\alpha}^{-1}x_{\alpha}(u)w_{\alpha})\\
=&e_{\alpha}(w_{\alpha})\sum_{u=\pm1}e_{\alpha}(x_{-\alpha}(u))\\
=&e_{\alpha}(w_{\alpha})(e_{\alpha}(w_{\alpha})+e_{\alpha}(w_{\alpha}^{-1}))\\
=&\epsilon\sqrt{2}e_{\alpha}(w_{\alpha}).
\end{align*}
The last equality follows from Lemma \ref{IrrationalWeyl}.\EndProof\\

\noindent \textbf{Remark:} The previous result implies that $(\frac{\epsilon}{\sqrt{2}}e_{\alpha}-2)(\frac{\epsilon}{\sqrt{2}}e_{\alpha}+1)=0$.

Now we will prove that $\mathcal{H}$ has an Iwahori-Matsumoto presentation. Let $H$ be the unital $\mathbb{C}$-algebra generated by the symbols $T_{w}$, where $w\in \widetilde{W}_{\mathrm{aff}}$, subject to the relations

\begin{align}
T_{w_{1}}T_{w_{2}}=&T_{w_{1}w_{2}}, \hspace{1cm} & \text{ for } w_{j}\in\widetilde{W}_{\mathrm{aff}}\text{ such that }\ell(w_{1}w_{2})=\ell(w_{1})+\ell(w_{2});\\
(T_{\alpha}+1)(T_{\alpha}-2)=&0,\hspace{1cm}&\text{ for any } \alpha\in \Delta\cup\{\alpha_{0}\}.
\end{align}
This is the algebra studied in Section 3 of \cite{L89}, where $v$ specialized to $\sqrt{2}$ and Lusztig's $L$ is equal to our $\ell$.

\begin{proposition}\label{AffineHeckePresentation}
The linear map $\Psi:H\rightarrow \mathcal{H}$ defined by $\Psi(T_{w})=(\frac{\epsilon}{\sqrt{2}})^{\ell(w)}e_{w}$ is an isomorphism of $\mathbb{C}$-algebras.
\end{proposition}

\textbf{Proof:} By Proposition \ref{LengthMult} and Proposition \ref{QuadRel} we know that $\Psi$ is a $\mathbb{C}$-algebra homomorphism. As $w$ varies over $\widetilde{W}_{\mathrm{aff}}$ the elements $T_{w}$ and $(\frac{\epsilon}{\sqrt{2}})^{\ell(w)}e_{w}$ form a basis for $H$ and $\mathcal{H}$ respectively. Thus $\Psi$ is an isomorphism of $\mathbb{C}$-algebras.\EndProof

\begin{proposition}\label{Invert}
For all $w\in \widetilde{W}_{\mathrm{aff}}$, the element $e_{w}\in \mathcal{H}$ is invertible. 
\end{proposition}

\textbf{Proof:} Let $w=\sigma w^{\prime}\in \widetilde{W}_{\mathrm{aff}}\cong \Omega\ltimes (W\ltimes 2Y)$, where $\sigma\in \Omega$ and $w^{\prime}\in W\ltimes 2Y$. The group $W\ltimes 2Y$ is an affine Weyl group with generators $w_{\alpha_{0}},\ldots,w_{\alpha_{r}}$. Suppose that $w^{\prime}=w_{\alpha_{i_{1}}}\ldots w_{\alpha_{i_{k}}}$ is a minimal expression for $w^{\prime}$ with respect to these generators. Then Proposition \ref{LengthMult} implies that $e_{w}=e_{\sigma}e_{\alpha_{i_{1}}}\ldots e_{\alpha_{i_{k}}}$. By Proposition \ref{QuadRel}, $e_{\alpha_{i_{j}}}$ is invertible. Since $e_{\sigma}e_{\sigma^{-1}}=1$, $e_{\sigma}$ is invertible. Thus $e_{w}$ is invertible.\EndProof

\textbf{Remark:}  It is this step where we deviate from Savin \cite{S04}. To establish the analog of Proposition \ref{Invert} in \cite{S04}, Savin studies the unramified principal series and invokes the theory of the Bernstein center. In the present paper, we can circumvent these issues using Theorem \ref{EqualLengths}.  Using this theorem we can prove an Iwahori-Matsumoto presentation, from which the invertibility is a simple consequence. In the tame case considered by Savin, the difficulty stems from the failure of the analog of equation (\ref{AffineDoubleCoset}). 

Let $\mathcal{H}_{0}$ be the subalgebra of $\mathcal{H}$ consisting of functions supported on $\widetilde{\textbf{G}}(\mathbb{Z}_{2})$. This subalgebra is generated by the elements $e_{\alpha}$, where $\alpha\in\Delta$. Let $\mathcal{A}$ be the subalgebra of $\mathcal{H}$ generated by $e_{\lambda}$, where $\lambda\in \widetilde{Y}$ and $\lambda$ is dominant.

\begin{proposition} \label{HeckeDecomp}
We have $\mathcal{H}=\mathcal{H}_{0}\cdot\mathcal{A}\cdot \mathcal{H}_{0}$.
\end{proposition}

\textbf{Proof:} The proof of Proposition 6.3 in Savin \cite{S04} directly adapts to this case.\EndProof

For $\lambda\in \widetilde{Y}$ let 
\[t_{\lambda}=q^{-\langle\lambda,\varrho\rangle}e_{\lambda_{1}}e_{\lambda_{2}}^{-1},\]
where $\lambda_1,\lambda_2\in \widetilde{Y}$ are dominant and $\lambda=\lambda_{1}-\lambda_{2}$. Note that the definition of $t_{\lambda}$ does not depend on the choice of $\lambda_{1}$ and $\lambda_{2}$.

Now we prove our main theorem. The Hecke algebra $\mathcal{H}$ has a Bernstein presentation.

\begin{theorem} \label{HeckeIso} 
Let $\mathscr{H}$ be the $\mathbb{C}$-algebra generated by $f_{\alpha}$, for all $\alpha\in \Delta$, and $u_{\lambda}$, for all $\lambda\in \widetilde{Y}$ modulo the relations 

\begin{enumerate}
\item $(f_{\alpha}-2)(f_{\alpha}+1)=0$;\label{QuadRelation}
\item \label{BraidRel} $\begin{cases}
f_{\alpha}\cdot f_{\beta}=f_{\beta}\cdot f_{\alpha},\hspace{4cm} \text{ if }\langle\alpha,\beta\rangle=0;\\
f_{\alpha}\cdot f_{\beta}\cdot f_{\alpha} =f_{\beta}\cdot f_{\alpha}\cdot f_{\beta},\hspace{3cm} \text{ if } \langle\alpha,\beta\rangle=-1;
\end{cases}$
\item \label{CoChar} $u_{\lambda}\cdot u_{\lambda^{\prime}}=u_{\lambda+\lambda^{\prime}}$.

\noindent Let $2m=\langle\alpha,\lambda\rangle$, where $\alpha$ is a simple root.
\item \label{HeckeAlgRelation}$\begin{cases}
f_{\alpha}\cdot u_{\lambda}=u_{\lambda}\cdot f_{\alpha}\hspace{5cm} \text{ if } m=0;\\
f_{\alpha}\cdot u_{\lambda}= u_{\lambda^{w_{\alpha}}}\cdot f_{\alpha}+\sum_{k=0}^{m-1}u_{\lambda-2k\alpha}, \hspace{2cm} \text{ if } m>0;\\
f_{\alpha}\cdot u_{\lambda}= u_{\lambda^{w_{\alpha}}}\cdot f_{\alpha}-\sum_{k=1}^{-m}u_{\lambda+2k\alpha}, \hspace{2cm} \text{ if } m<0.

\end{cases}$
\end{enumerate}

The map $A:\mathscr{H}\rightarrow \mathcal{H}$ defined by 
\[\begin{cases}
A(u_{\lambda})=t_{\lambda}\\
A(f_{\alpha})=\frac{\epsilon}{\sqrt{2}}e_{w_{\alpha}}
\end{cases}\]
is an isomorphism of $\mathbb{C}$-algebras.
\end{theorem}

\textbf{Proof:} For the Hecke algebra $\mathcal{H}$, we established the analog identity (\ref{QuadRelation}), in Proposition \ref{QuadRel},  and the analogs of identities (\ref{BraidRel}) and (\ref{CoChar}), in Proposition \ref{LengthMult}. Identity (\ref{HeckeAlgRelation}) in $\mathcal{H}$ is exactly Proposition 3.6 in \cite{L89}, which is applicable by Proposition \ref{AffineHeckePresentation}. 

We have shown that $\mathcal{H}$ satisfies identities (\ref{QuadRelation}), (\ref{BraidRel}), (\ref{CoChar}), and $(\ref{HeckeAlgRelation})$. Thus the map $A:\mathscr{H}\rightarrow\mathcal{H}$ is a well-defined homomorphism of algebras. Proposition \ref{HeckeDecomp} implies that $A$ is surjective, and to prove that $A$ is injective one can apply the trick of Lemma 7.6 and the discussion immediately follow it in Savin \cite{S04}. Thus the map $A:\mathscr{H}\rightarrow\mathcal{H}$ is is an isomorphism of algebras.\EndProof

\begin{corollary}\label{tensorDecomp} The subalgebra generated by all of the $t_{\lambda}$ is isomorphic to the group algebra $\mathbb{C}[\widetilde{Y}]$ and the multiplication map defines an isomorphism of vector spaces
\begin{equation}
\mathbb{C}[\widetilde{Y}]\otimes \mathcal{H}_{0}\cong \mathcal{H}.
\end{equation}

\end{corollary}


\section{Representations Generated By $\tau$-isotypic Vectors}\label{BernsteinComponent}

Let $\mathcal{R}(\widetilde{G})$ be the category of smooth $\widetilde{G}$-representations, and let $\mathcal{R}(\widetilde{G},\tau)$ be the full subcategory of representations that are generated by their $\tau$-isotypic vectors.

Now we prove that the Bernstein component of the genuine unramified principal series is $\mathcal{R}(\widetilde{G},\tau)$.

\begin{theorem}\label{BernsteinBlock}$ $
\begin{enumerate}
\item\label{BB1} The category $\mathcal{R}(\widetilde{G},\tau)$ is closed relative to subquotients in $\mathcal{R}(\widetilde{G})$. 

\item\label{BB2} For every irreducible object $(\pi,\mathscr{V})$ in $\mathcal{R}(\widetilde{G},\tau)$, there is an unramified character $\chi:T\rightarrow \mathbb{C}^{\times}$ such that $(\pi,\mathscr{V})$ is isomorphic to an irreducible subquotient of $\mathrm{Ind}_{\widetilde{B}}^{\widetilde{G}}(i(\chi))$. 
\end{enumerate}
\end{theorem}

We begin with a few preliminary results. Recall that $\tau^{\diamond}$ is a genuine irreducible representation of $\widetilde{T}$, and $\tau$ is the inflation of $\tau^{\diamond}$ to $\Gamma_{0}(4)$. We will write $\tau_{\widetilde{T}}$ for the restriction $\tau$ to $\widetilde{T}\cap \widetilde{\Gamma}_{0}(4)$.

\begin{lemma}\label{CompactTType}
Let $\pi$ be a smooth irreducible $\widetilde{T}$-representation. Then the $\tau_{\widetilde{T}}$-isotypic subspace $\pi^{\tau_{\widetilde{T}}}\neq0$ if and only if $\pi\cong i(\chi)$ as $\widetilde{T}$-representations, for some unramified character $\chi:T\rightarrow \mathbb{C}^{\times}$.
\end{lemma}

\textbf{Proof:} Suppose that $\pi^{\tau_{\widetilde{T}}}\neq0$. Since $\widetilde{T}\cap \widetilde{\Gamma}_{0}(4)\cong \widetilde{T}^{\diamond}\times T_{1}^{*}$ and $\tau_{\widetilde{T}}\cong \tau^{\diamond}\otimes 1$ under this isomorphism, the central character of $\pi$ is equal to the central character of $i(\chi)$ for some unramified $\chi$. Since an irreducible representation of $\widetilde{T}$ in which $T_{1,8}^{*}$ acts by the identity is determined by its central character, we have $\pi\cong i(\chi)$.

Conversely, suppose that $\pi\cong i(\chi)$. Consider the subspace $\pi^{T_{1}^{*}}$ of $T_{1}^{*}$-invariants. It is a $\widetilde{T}^{\diamond}T_{1}^{*}$-subrepresentation because $\widetilde{T}\cap\widetilde{\Gamma}_{0}(4)\cong \widetilde{T}^{\diamond}\times T_{1}^{*}$. We can see that $\pi^{T_{1}^{*}}\neq 0$ by showing that $V(\gamma_{2})^{T_{1}^{*}}\neq 0$. This can be seen by viewing $V(\gamma_{2})$ as an induced representation. Specifically, let $\gamma_{2}^{\prime}$ be the character of the maximal abelian subgroup $T_{1}^{*}\Upsilon(\widetilde{Y})\mu_{2}$ extending $\gamma_{2}$ such that $\gamma_{2}^{\prime}|_{T_{1}^{*}}=1$. Then $V(\gamma_{2})\cong\mathrm{Ind}_{T_{1}^{*}\Upsilon(\widetilde{Y})\mu_{2}}^{T_{1}^{*}\Upsilon(Y)\mu_{2}}(\gamma_{2}^{\prime})$ and the functions supported on $T_{1}^{*}\Upsilon(\widetilde{Y})\mu_{2}$ are fixed by $T_{1}^{*}$.\EndProof

\begin{proposition}\label{JacquetIso}
Let $(\pi,\mathscr{V})$ be a smooth $\widetilde{G}$-representation. We will write $(\pi_{U},\mathscr{V}_{U})$ for the normalized Jacquet module with respect to $U^{*}$, and $q:\mathscr{V}\rightarrow \mathscr{V}_{U}$ for the canonical quotient map. Then the induced map
\begin{equation}
q:\mathscr{V}^{\tau}\rightarrow \mathscr{V}_{U}^{\tau_{\widetilde{T}}}
\end{equation}
is an isomorphism.
\end{proposition}

\textbf{Proof:} This follows from Theorem 7.9 in Bushnell-Kutzko \cite{BK98}. Hypothesis a) in Theorem 7.9 holds because $\widetilde{\Gamma}_{0}(4)$ possess an Iwahori-factorization; hypothesis b) holds because of our Proposition \ref{Invert}.\EndProof

\begin{proposition}\label{Gamma04Type}
Let $\pi$ be a smooth irreducible $\widetilde{G}$-representation. Then the $\tau$-isotypic space $\pi^{\tau}\neq 0$ if and only if $\mathrm{Hom}_{\widetilde{G}}(\pi,\mathrm{Ind}_{\widetilde{B}}^{\widetilde{G}}i(\chi))\neq0$ for some unramified character $\chi:T\rightarrow \mathbb{C}^{\times}$.
\end{proposition}

\textbf{Proof:} Suppose that the $\pi^{\tau}\neq 0$. By Proposition \ref{JacquetIso}, $\pi_{U}^{\tau_{\widetilde{T}}}\cong \pi^{\tau}\neq 0$. Since $\pi$ is irreducible, $\pi_{U}$ is finitely generated. Because $\pi_{U}^{\tau_{\widetilde{T}}}\neq 0$ we can apply Lemma \ref{CompactTType} and the Bernstein decomposition to see that $\pi_{U}$ has an irreducible quotient isomorphic to $i(\chi)$ for some unramified character $\chi:T\rightarrow \mathbb{C}^{\times}$. Thus by Frobenius reciprocity $0\neq \mathrm{Hom}_{\widetilde{T}}(\pi_{U},i(\chi))\cong \mathrm{Hom}_{\widetilde{G}}(\pi,\mathrm{Ind}_{\widetilde{B}}^{\widetilde{G}}i(\chi))$.


Conversely, suppose that $\mathrm{Hom}_{\widetilde{G}}(\pi,\mathrm{Ind}_{\widetilde{B}}^{\widetilde{G}}i(\chi))\neq0$. By Frobenius reciprocity there is a surjective map $\psi:\pi_{U}\rightarrow i(\chi)$ of $\widetilde{T}$-representations, since $i(\chi)$ is irreducible. By Lemma \ref{CompactTType}, $\tau_{\widetilde{T}}$ is a $\widetilde{T}\cap\widetilde{\Gamma}_{0}(4)$-subrepresentation of $i(\chi)$. Since $\tau_{\widetilde{T}}$ is an irreducible representation of the compact group $\widetilde{T}\cap\widetilde{\Gamma}_{0}(4)$, $\tau_{\widetilde{T}}$ must be a $\widetilde{T}\cap \widetilde{\Gamma}_{0}(4)$-subrepresentation of $\pi_{U}$. Thus $\pi_{U}^{\tau_{\widetilde{T}}}\neq0$. Finally, by Proposition \ref{JacquetIso}, $\pi^{\tau}\cong \pi_{U}^{\tau_{\widetilde{T}}}\neq 0$.

\textbf{Proof of Theorem \ref{BernsteinBlock}:} We begin with item (\ref{BB1}). Let $(\pi,\mathscr{V})$ be smooth $\widetilde{G}$-representation generated by $\mathscr{V}^{\tau}$. Using the Bernstein decomposition of $\mathcal{R}(\widetilde{G})$ it suffices to assume that $(\pi,\mathscr{V})$ is contained in a single Bernstein component. Under this assumption, we claim that $\mathscr{V}$ is in the block associated to the $\widetilde{T}$-representation $i(1)$. To see this we consider $\mathscr{W}$ a $\widetilde{G}$-representation generated by some $v\in \mathscr{V}^{\tau}-\{0\}$. Since $\mathscr{W}$ is finitely generated it has an irreducible quotient $\mathscr{W}^{\prime}$. Since the image of $v$ is nonzero in $\mathscr{W}^{\prime}$, it follows that $(\mathscr{W}^{\prime})^{\tau}\neq 0$. Thus by Proposition \ref{Gamma04Type}, $\mathscr{V}$ has an irreducible subquotient in the block associated to the $\widetilde{T}$-representation $i(1)$, and it follows that $\mathscr{V}$ must be in the same block.

Now suppose that $\mathscr{W}$ is any subrepresentation of $\mathscr{V}$ such that $\mathscr{W}^{\tau}$ does not generate $\mathscr{W}$. It suffices to assume that $\mathscr{W}^{\tau}=0$, since we can take the quotient of $\mathscr{V}$ and $\mathscr{W}$ by the $\widetilde{G}$-subrepresentation generated by $\mathscr{W}^{\tau}$. Under this assumption, we claim that $\mathscr{W}$ has an irreducible subquotient that does not live in the block associated to the $\widetilde{T}$-representation $i(1)$. This follows from Proposition \ref{Gamma04Type}, since the $\tau$-isotypic subspace of any irreducible quotient of $\mathscr{W}$ must be trivial. This implies that $\mathscr{W}$ and $\mathscr{V}$ are in different blocks. This contradiction prove item (\ref{BB1}). 

Item (\ref{BB2}) is exactly Proposition \ref{Gamma04Type}.\EndProof

\textbf{Remark:} In the language of Bushnell-Kutzko \cite{BK98}, Theorem \ref{BernsteinBlock} states that $\mathcal{R}(\widetilde{G},\tau)$ is the Bernstein component associated to the inertial support $[\widetilde{T},i(1)]_{\widetilde{G}}$; Proposition \ref{Gamma04Type} states that $(\widetilde{\Gamma}_{0}(4),\tau)$ is a $[\widetilde{T},i(1)]_{\widetilde{G}}$-type. Moreover, Theorem 4.3 in Bushnell-Kutzko \cite{BK98} states that the functor $V\mapsto \mathrm{Hom}_{\widetilde{\Gamma}_{0}(4)}(\tau,V)$ defines an equivalence of categories between $\mathcal{R}(\widetilde{G},\tau)$ and $\mathcal{H}$-mod. Thus the Bernstein component associated to $[\widetilde{T},i(1)]_{\widetilde{G}}$ is equivalent to $\mathcal{H}$-mod.


\section{Local Shimura Correspondence}\label{ShimuraCor}

Consider the split algebraic group $G^{\prime}=G/Z_{2}$, where $Z_{2}$ is the $2$-torsion of the center of $G$. The possible $2$-groups which can arise as $Z_{2}$ are described in Table \ref{Central2Tor}. Let $I^{\prime}$ denote an Iwahori subgroup of $G^{\prime}$. The Bernstein Presentation of the Iwahori-Hecke algebra $\mathscr{H}(G^{\prime},I^{\prime})$ of  $G^{\prime}$ implies the vector space isomorphism $\mathcal{H}(G^{\prime},I^{\prime})\cong \mathbb{C}[Y^{*}\cap \frac{1}{2}Y]\otimes H$, since $Y^{*}\cap \frac{1}{2}Y$ is the co-weight lattice of $G^{\prime}$. However, $2(Y^{*}\cap \frac{1}{2}Y)=\widetilde{Y}$, so $\lambda\mapsto2\lambda$ defines an algebra isomorphism $\mathbb{C}[Y^{*}\cap \frac{1}{2}Y]\rightarrow \mathbb{C}[\widetilde{Y}]$, which extends to an algebra isomorphism $\mathscr{H}(G^{\prime},I^{\prime})\stackrel{\phi}{\rightarrow}\mathcal{H}$, by Theorem \ref{HeckeIso} and Corollary \ref{tensorDecomp}.

\begin{center}
\begin{tabular}{|c|c|c|c|c|c|c|c|}\hline 
$\Phi$ & $A_{2n}$ & $A_{2n+1}$ & $D_{2n}$ & $D_{2n+1}$ & $E_{6}$ & $E_{7}$ & $E_{8}$\\ \hline
$Z_{2}$ & $\{1\}$ & $\mathbb{Z}/2\mathbb{Z}$ & $\mathbb{Z}/2\mathbb{Z}^{\oplus2}$ & $\mathbb{Z}/2\mathbb{Z}$ & $\{1\}$ & $\mathbb{Z}/2\mathbb{Z}$ & $\{1\}$\\ \hline
\end{tabular}
\end{center}
 \begingroup
 \captionof{table}{Central $2$-torsion of simple, simply-laced, simply-connected Chevalley groups.} \label{Central2Tor}
 \endgroup

Now we can describe the Shimura correspondence. Let $\mathcal{R}(G^{\prime},I^{\prime})$ be the category of smooth representations of $G^{\prime}=G/Z_{2}$ that are generated by their $I^{\prime}$ fixed vectors. Recall that $\mathcal{R}(\widetilde{G},\tau)$ is the category of smooth representations of $\widetilde{G}$ that are generated by their $\tau$-isotypic vectors.

\begin{theorem}\label{main}
The isomorphism $\mathscr{H}(G^{\prime},I^{\prime})\stackrel{\phi}{\rightarrow}\mathcal{H}$ induces an equivalence of categories between $\mathcal{R}(G^{\prime},I^{\prime})$ and $\mathcal{R}(\widetilde{G},\tau)$.
\end{theorem}

\textbf{Proof:} The algebra isomorphism $\mathscr{H}(G^{\prime},I^{\prime})\stackrel{\phi}{\rightarrow}\mathcal{H}$ induces an equivalence of categories between $\mathscr{H}(G^{\prime},I^{\prime})$-mod and $\mathcal{H}$-mod.

By results of Bushnell-Kutzko \cite{BK98} (Theorem 4.3), the functor $V\mapsto \mathrm{Hom}_{I^{\prime}}(1,V)=V^{I^{\prime}}$ defines an equivalence of categories between $\mathcal{R}(G^{\prime},I^{\prime})$ and $\mathscr{H}(G^{\prime},I^{\prime})$-mod. Thus, $\mathcal{R}(G^{\prime},I^{\prime})$ is equivalent to $\mathcal{H}$-mod.

Again, by results of Bushnell-Kutzko \cite{BK98} (Theorem 4.3), the functor $V\mapsto \mathrm{Hom}_{\widetilde{\Gamma}_{0}(4)}(\tau,V)$ defines an equivalence of categories between $\mathcal{R}(\widetilde{G},\tau)$ and $\mathcal{H}$-mod. Theorem 4.3 is applicable because our Proposition \ref{Gamma04Type} implies that $(\widetilde{\Gamma}_{0}(4),\tau)$ is a $[\widetilde{T},i(1)]_{\widetilde{G}}$-type. Thus $\mathcal{R}(G^{\prime},I^{\prime})$ is equivalent to $\mathcal{R}(\widetilde{G},\tau)$.\EndProof


\textbf{Remark:} Finally, we enumerate the choices we made in constructing $\mathscr{H}(G^{\prime},I^{\prime})\stackrel{\phi}{\rightarrow}\mathcal{H}$. Ostensibly we have made two choices. However, only one of these choices is material, the choice of a genuine character of $Z(\widetilde{T}^{\diamond})\cong \mathrm{Hom}(\widetilde{Y}/2Y,\mathbb{C}^{\times})$. 

First, the Hecke algebra $\mathcal{H}$ depends on the choice of an irreducible genuine Weyl group invariant $\widetilde{T}^{\diamond}$-representation $(\tau^{\diamond},E)$. Since each irreducible genuine $\widetilde{T}^{\diamond}$-representation is Weyl group invariant and determined by its central character, we can count the number of genuine central characters of $\widetilde{T}^{\diamond}$. Note that 
\begin{equation}
1\rightarrow \mu_{2}\rightarrow \widetilde{T}^{\diamond}\rightarrow Y\otimes \mu_{2}\cong Y/2Y\rightarrow 1.
\end{equation}
One can show that the center of $\widetilde{T}^{\diamond}$ is the preimage of $\widetilde{Y}/2Y$. So we get the split exact sequence of $\mathbb{F}_{2}$-vector spaces
\begin{equation}
1\rightarrow \mu_{2}\rightarrow Z(\widetilde{T}^{\diamond})\rightarrow  \widetilde{Y}/2Y\rightarrow 1.
\end{equation}
Thus we see that genuine characters of $Z(\widetilde{T}^{\diamond})$ are in bijection with Hom$(\widetilde{Y}/2Y,\mathbb{C}^{\times})$. An explicit description of $\widetilde{Y}/2Y$ can be found in Section 16.1 in Gan-Gao \cite{GG18}. We note that $Z_{2}\cong\widetilde{Y}/2Y$.

Second, we chose a particular normalization for the basis elements $e_{w}\in\mathcal{H}$. (Recall the discussion before Proposition \ref{BasisWellDefined}.) This normalization depends on a choice of an extension of $\tau^{\diamond}$ to $\mathcal{W}$, which we called $\tau_{\mathcal{W}}$. The number of such extensions is in bijection with the set $\mathrm{Hom}(\mathcal{W}/\widetilde{T}^{\diamond},\mathbb{C}^{\times})\cong\mathrm{Hom}(W,\mathbb{C}^{\times})$. We claim that $\mathrm{Hom}(W,\mathbb{C}^{\times})\cong \mathbb{Z}/2\mathbb{Z}$, where the nontrivial character is defined by $w\stackrel{\mathrm{sign}}{\mapsto} (-1)^{\ell(w)}$. This is the result of the following standard facts about Weyl groups: the group $W$ is generated by simple reflections; the Weyl group acts transitively on roots of the same length; for any $\alpha,\beta\in\Phi$, $w_{\alpha}w_{\beta}w_{\alpha}^{-1}=w_{w_{\alpha}(\beta)}$. Furthermore, we see that if $\tau_{\mathcal{W}}$ is one of the extensions, then the other must be given by $\tau_{\mathcal{W}}\otimes \mathrm{sign}$, where now we write $\mathrm{sign}$ for the inflation of the sign representation of $W$ to $\mathcal{W}$. If we replace $\tau_{\mathcal{W}}$ by $\tau_{\mathcal{W}}\otimes \mathrm{sign}$, this has the effect of changing $\epsilon$ to $-\epsilon$; $e_{\alpha}$ to $-e_{\alpha}$, where $\alpha\in\Delta$; and $t_{\lambda}$ remains invariant for $\lambda\in \widetilde{Y}$. Thus $\frac{\epsilon}{\sqrt{2}}e_{\alpha}$ and $t_{\lambda}$ remain invariant for $\alpha\in \Delta$ and $\lambda\in \widetilde{Y}$. This implies that the isomorphism of Theorem \ref{HeckeIso}, and thus the isomorphism $\mathscr{H}(G^{\prime},I^{\prime})\stackrel{\phi}{\rightarrow}\mathcal{H}$, does not depend on our choice of extension of $\tau^{\diamond}$ from $\widetilde{T}^{\diamond}$ to $\mathcal{W}$.

In summary, the isomorphism $\mathscr{H}(G^{\prime},I^{\prime})\stackrel{\phi}{\rightarrow}\mathcal{H}$ only depends on the choice of a genuine character of $Z(\widetilde{T}^{\diamond})\cong \mathrm{Hom}(\widetilde{Y}/2Y,\mathbb{C}^{\times})$.









\begin{thebibliography}{9}

\bibitem{ABPTV07}
Adams, J.; Barbasch, D.; Paul, A.; Trapa, P.; Vogan, D. A., Jr. 
\textit{Unitary Shimura correspondences for split real groups}. 
J. Amer. Math. Soc. 20 (2007), no. 3, 701--751.

\bibitem{BD84}
Bernstein, J.; Deligne, P.
\textit{Le ``Centre'' de Bernstein.}
Repr\'{e}sentations des groupes r\'{e}ductifs sur un corps local. Hermann, Paris 1984.

\bibitem{BK93}
Bushnell, Colin J.; Kutzko, Philip C. 
\textit{The admissible dual of ${\rm GL}(N)$ via compact open subgroups.} 
Annals of Mathematics Studies, 129. Princeton University Press, Princeton, NJ, 1993. {\rm xii}+313 pp. ISBN: 0-691-03256-4; 0-691-02114-7

\bibitem{BK98}
Bushnell, Colin J.; Kutzko, Philip C. 
\textit{Smooth representations of reductive $p$-adic groups: structure theory via types.} 
Proc. London Math. Soc. (3) 77 (1998), no. 3, 582--634.

\bibitem{GG18}
Gan, Wee Teck; Gao, Fan. 
\textit{The Langlands-Weissman program for Brylinski-Deligne extensions. L-groups and the Langlands program for covering groups.}
 Astérisque 2018, no. 398, 187--275. ISBN: 978-2-85629-845-9

\bibitem{H75}
Humphreys, James E. 
\textit{Linear algebraic groups.} 
Graduate Texts in Mathematics, No. 21. Springer-Verlag, New York-Heidelberg, 1975. {\rm xiv}+247 pp.

\bibitem{IM65}
Iwahori, N.; Matsumoto, H. 
\textit{On some Bruhat decomposition and the structure of the Hecke rings of $p$-adic Chevalley groups}. 
Inst. Hautes Études Sci. Publ. Math. No. 25 1965 5--48.

\bibitem{LS10a}
Loke, Hung Yean; Savin, Gordan. 
\textit{Modular forms on non-linear double covers of algebraic groups}. 
Trans. Amer. Math. Soc. 362 (2010), no. 9, 4901--4920. 

\bibitem{LS10b}
Loke, Hung Yean; Savin, Gordan. 
\textit{Representations of the two-fold central extension of ${\rm SL}_2(\Bbb Q_2)$}. 
Pacific J. Math. 247 (2010), no. 2, 435--454. 

\bibitem{L89}
Lusztig, George. 
\textit{Affine Hecke algebras and their graded version.}
 J. Amer. Math. Soc. 2 (1989), no. 3, 599--635.

\bibitem{M69}
Matsumoto, Hideya. 
\textit{Sur les sous-groupes arithmétiques des groupes semi-simples déployés.} (French) 
Ann. Sci. École Norm. Sup. (4) 2 (1969), 1--62.

\bibitem{M68}
Moore, Calvin C. 
\textit{Group extensions of $p$-adic and adelic linear groups.} 
Inst. Hautes Études Sci. Publ. Math. No. 35, 1968 157--222.

\bibitem{S88}
Savin, Gordan. 
\textit{Local Shimura correspondence}. 
Math. Ann. 280 (1988), no. 2, 185--190.

\bibitem{S04}
Savin, Gordan. 
\textit{On unramified representations of covering groups}. 
J. Reine Angew. Math. 566 (2004), 111--134.

\bibitem{S73}
Stein, Michael R. 
\textit{Surjective stability in dimension $0$ for $K\sb{2}$ and related functors}. 
Trans. Amer. Math. Soc. 178 (1973), 165--191. 

\bibitem{S68}
Steinberg, Robert. 
\textit{Lectures on Chevalley groups}. Notes prepared by John Faulkner and Robert Wilson. Revised and corrected edition of the 1968 original [MR0466335]. With a foreword by Robert R. Snapp. 
University Lecture Series, 66. American Mathematical Society, Providence, RI, 2016. xi+160 pp. ISBN: 978-1-4704-3105-1

\bibitem{TW18}
Takeda, Shuichiro; Wood, Aaron. 
\textit{Hecke algebra correspondences for the metaplectic group.} 
Trans. Amer. Math. Soc. 370 (2018), no. 2, 1101--1121.

\bibitem{W14}
Wood, Aaron. 
\textit{A minimal even type of the 2-adic Weil representation.} 
Math. Z. 277 (2014), no. 1-2, 257--283.







\end{thebibliography}
\end{document}